\documentclass[epsf]{article}
\usepackage{amsmath,amsthm,amsfonts,latexsym,amscd}\usepackage{amsmath,amsthm,amsfonts,latexsym,amscd}
\title{Submanifolds of Constant Negative Curvature:\\A Generalization of Hilbert\rq s Theorem}
\author{John Douglas Moore\\Department of 
Mathematics\\University of California\\
Santa Barbara, CA, USA 93106\\e-mail: dmoore@ucsb.edu}
\date{}

\begin{document}

\maketitle

\begin{abstract}
We use the generalized Gauss-Bonnet formula for Riemannian polyhedra discovered by Allendoerfer, Weil and Chern to show that hyperbolic space of dimension $n$ has no isometric immersion into Euclidean space of dimension $2n-1$.
\end{abstract}

\section{Introduction}

In the nineteenth century, two noneuclidean geometries were discovered, spherical and hyperbolic.  The spherical geometry has an obvious model as the geometry of the standard unit sphere in three-dimensional Euclidean space ${\mathbb E}^3$, but hyperbolic geometry was more mysterious.  At the beginning of the twentieth century, David Hilbert \cite{Hi} proved that the hyperbolic plane cannot be realized as the Riemannian geometry induced on an immersed surface in ${\mathbb E}^3$.

There are two steps to the standard proof of Hilbert\rq s theorem, as presented at the very end of do Carmo's book (\S 5.11 of \cite{dC}).  The first step shows that if there exists an immersed surface in ${\mathbb E}^3$ isometric to the hyperbolic plane, then it must be covered by a global coordinate system $(z,w)$ whose coordinate vectors are unit-length asymptotic vectors.  (This is what is sometimes called a Chebyshev net, and is related to the design of clothing as described in \cite{G}.)  When expressed in these coordinates, the condition that the Gaussian curvature be $-1$ implies via a short calculation that the angle $\theta $ between the curves tangent to the asymptotic vectors satisfies the {\em sine-Gordon equation\/}
$$\frac{\partial ^2 \theta}{\partial z \partial w} = \sin \theta,$$
an angle which must satisfy the inequalities $0 < \theta < \pi$.  In the second step, one shows that the sine-Gordon equation has no global solutions which satisfy the inequalities $0 < \theta < \pi$, either by an elementary argument due to Holmgren (\cite{Ho} or \cite{McC}, pages 231-232) or by means of a formula discovered by the Greek mathematician Ioannis Hazzidakis \cite{Haz}.  Indeed, one can verify that the integral of $\sin \theta $ over any coordinate rectangle, defined in terms of the asymptotic coordinates by
\begin{equation} R = \{ (z,w) : a \leq z \leq b, \quad c \leq w \leq d \}, \label{E:asymptoticrectangle} \end{equation}
is the hyperbolic area of that rectangle, and then use the two-dimensional Gauss-Bonnet formula to show that the sine-Gordon equation is equivalent to the {\em Hazzidakis formula\/}
\begin{equation} \hbox{Hyperbolic area of $R$} = \theta (b,d) - \theta (b,c) - \theta (a,d) + \theta (a,c), \label{E:hilb11} \end{equation}
for all such coordinate rectangles.  This shows that the area of the coordinate rectangle $R$ is bounded by $2\pi $, contradicting the fact that since the hyperbolic plane has infinite area, coordinate rectangles within it must have arbitrarily large area.  This finishes our sketch of the proof of Hilbert\rq s original theorem.  Note that we can rewrite (\ref{E:hilb11}) as
\begin{equation} \frac{1}{2\pi} \left[ \hbox{Hyperbolic area of $R$} \right] = \sum_{i=1}^4 \frac{\varepsilon _i}{2\pi} - 1, \label{E:hilb12} \end{equation}
where the $\varepsilon _i$'s are the four interior (or exterior) angles at the vertices of the coordinate rectangle.  We will recall a Gauss-Bonnet argument that establishes (\ref{E:hilb12}) at the end of \S \ref{S:integrationoverfibers}.

We emphasize that the argument via the Hazzidakis formula actually proves a quantitative refinement of Hilbert\rq s theorem:  A surface in ${\mathbb E}^3$ of constant negative curvature cannot contain an asymptotic coordinate rectangle of area larger than $2\pi $.

Hilbert's theorem contrasts with the fact that noncomplete surfaces of constant negative curvature are plentiful.  Indeed, it can be shown that they correspond to local solutions of the sine-Gordon equation, and using the Cauchy-Kovalevskaya Theorem, one can show that real analytic solutions depend upon two functions of a single variable (just like the linear wave equation in which we replace $\sin \theta $ by zero).

The purpose of this article is to explain how Hilbert's theorem extends to hyperbolic spaces of higher dimension.

Here is some background on this problem.  \`Elie Cartan \cite{Ca} proved that there are no isometric immersions from open subsets of the $n$-dimensional hyperbolic space ${\mathbb H}^n$ into ${\mathbb E}^N$ when $N < 2n-1$, and that solutions to the local problem of finding such isometric immersions depend upon $n(n-1)$ functions of a single variable when $N = 2n-1$.  He established this as a consequence of Cartan-K\"ahler theory for local real analytic solutions to nonlinear systems of partial differential equations, which is presented in \cite{Ca2} and \cite{BCGGG}.  Thus there is a rich theory of such local solutions, and it bears a striking resemblance to the theory of surfaces of constant negative curvature in ${\mathbb E}^3$.  Such solutions can be regarded as solutions to a generalization of the sine-Gordon equation, and just as for surfaces of constant negative curvature, new examples can be generated from old by B\"acklund transformations.  These facts are explained in articles by Tenenblat and Terng \cite{TT} and Terng \cite{T}.

It was conjectured in \cite{Mo69} and \cite{Mo} that there is no global isometric immersion from ${\mathbb H}^n$ into ${\mathbb E}^{2n-1}$, and the first step of a projected proof was established there: the existence of a global Chebyshev net generated by unit--length asymptotic vector fields.  This conjecture was mentioned in one of Yau's problem lists; see \cite{Y}, page 682.  Positive evidence for the conjecture was provided by Xavier \cite{X}, who showed that ${\mathbb H}^n/\Gamma$ cannot be isometrically immersed in ${\mathbb E}^{2n-1}$, when $\Gamma $ is a non-elementary group of isometries, and by a later result of Nikolayevsky \cite{N} for the case in which $\Gamma $ is nontrivial.  We can refer to \cite{FX} and \cite{DOV} for additional related recent literature.

In this article, we will prove the conjecture stated above:

\vskip .1 in
\noindent
{\bf Main Theorem.} {\sl If $n \geq 2$, there is no isometric immersion from ${\mathbb H}^{n}$ into ${\mathbb E}^{2n-1}$.}

\vskip .1 in
\noindent
To prove this when $n = 2m$ is even, we establish a generalization of the Hazzidakis formula to isometric immersions from an open subset $U$ of ${\mathbb H}^n$ into ${\mathbb E}^{2n-1}$ in this case, as a consequence of the generalized Gauss-Bonnet formula for Riemannian polyhedra developed by Allenndoerfer and Weil \cite{AW}.   We will apply the Allenndoerfer-Weil formula to a suitable family of \lq\lq asymptotic polyhedra" $P_a$ within $U$ defined for $a \in (0, \infty)$, which specialize to coordinate rectangles in the case of surfaces.  Each such $P_a$ will be homeomorphic to an $n$-disk $D^n$ in ${\mathbb R}^n$ with a boundary $\partial P_a$ that is a simplicial complex homeomorphic to $S^{n-1}$; thus $P_a$ is an imbedded $n$-cell in the terminology utilized in \S 6 of Allendoerfer and Weil \cite{AW}.  However, unlike the Hazzidakis formula for surfaces in ${\mathbb E}^3$, the $n$-dimensional formula obtained from the Allendoerfer-Weil formula applied to the Levi-Civita connection contains curvature terms of every even dimension between zero and $n$ which our methods do not show to be zero.  Instead of attempting to evaluate these terms, we average the Allendoerfer-Weil formulae over various metric connections which are related to the Levi-Civita connection by the action of the symmetry group $G$ of $P_a$.  We show that the averaged Allendoerfer-Weil formula only contains nonzero terms in dimensions zero and $n$, thus yielding a higher-dimensional version of the Hazzidakis formula.  

Just as in the case of surfaces, the family of asymptotic polyhedra exhaust ${\mathbb H}^n$ as $a \rightarrow \infty$, and this contradicts the infinite volume of ${\mathbb H}^n$, thus providing a proof of the Main Theorem in the even-dimensional case.  As in the case of surfaces, the Hazzidakis formula gives a quantitative refinement: a bound on the volume of an asymptotic polyhedron contained in a $n$-dimensional submanifold of constant negative curvature in ${\mathbb E}^{2n-1}$.

This Hazzidakis argument doesn't work when $n$ is odd, since the Gauss-Bonnet integrand in the $n$-dimensional term of the Allendoerfer-Weil formula is zero for odd-dimensional Riemannian manifolds.  However, in this case we can apply the Gauss-Bonnet formula for manifolds with boundary, and we find that it is possible to develop a similar formula for the $(n-1)$-dimensional volume of the boundary $\partial P_a$ of $P_a$ in terms of lower dimensional integrals, the boundary being even-dimensional in this case.  As in the even-dimensional case, we need to average over connections obtained from the Levi-Civita connection by the $G$-action to eliminate terms of nonzero dimension lower than $n-1$ and obtain a bound on the volume of $\partial P_a$.  We can then apply an isoperimetric inequality (Theorem 34.2.6 in Burago and Zallgaller \cite{BZ}) to achieve an upper bound on the volume enclosed within $P_a$, contradicting once again the fact that as $a \rightarrow \infty$, the polyhedra $P_a$ exhaust the infinite volume hyperbolic space ${\mathbb H}^n$, thereby establishing the Main Theorem in the odd-dimensional case.

\vskip .1 in
\noindent
{\bf Outline of the remainder of this article.}  In \S \ref{S:adaptedframe}, we present the local theory of $n$-dimensional submanifolds of constant negative curvature in ${\mathbb E}^{2n-1}$, including a proof of the existence of Chebyshev nets (improving slightly our earlier arguments in \cite{Mo}).  In \S \ref{S:asymptotic}, we describe the asymptotic polyhedra $P_a$ in more detail, along with the symmetry group $G$ which acts on these polyhedra, and describe the family of metric connections in the tangent bundle that can be obtained from the Levi-Civita connection by the action of $G$.

Since the Gauss-Bonnet formula is crucial to our arguments, we review Chern's intrinsic proof of this formula in some detail in \S \ref{S:gaussbonnet}, using an argument based on construction of a Thom form that is adapted from Matthai and Quillen \cite{MQ}.  In \S\ref{S:integrationoverfibers} we extend the argument to the Gauss-Bonnet formula for manifolds with smooth boundary by integrating this Thom form over fibers, using constructions presented in Bott \cite{B}.  This approach is slightly more conceptual than the approach of Allendoerfer, Weil and Chern which used Weyl's calculation of the volume of tubes (as explained in \cite{W} and \cite{Gr}).  The argument we present provides a geodesic curvature form that depends only on a unit-length vector field homotopic to the outward-pointing unit normal along the boundary, which is not necessarily normal to the boundary, a fact essential to our proof.\footnote{Since integration along fibers in this case is integration along the rays generated by the vector field, the geodesic curvature form could also be expressed in terms of a {\em nowhere zero\/} vector field along the boundary homotopic to the outward-pointing unit normal.}

To obtain a Hazzidakis formula for $P_a$, we need to extend this Gauss-Bonnet formula to Riemannian polyhedral $n$-cells as described by Allendoerfer and Weil, and we do this by further integration over the fiber.  In \S \ref{S:dualCW}, we describe dual CW decompositions of the boundary $\partial P_a$ of the asymptotic polyhedron $P_a$  and the $(n-1)$-sphere fiber of the unit tangent bundle, and examine the action of the group $G$ of symmetries on the cellular decomposition of the fiber of the unit tangent bundle.  We then establish the Allendoerfer-Weil formula for our specific asymptotic polyhedra in \S \ref{S:faces}, noting that it applies to any metric connection on the tangent bundle to $M$.  Finally, we state the Hazzidakis formulae for both $P_a$ and $\partial P_a$ in \S \ref{S:hazzidakis}, and note that establishing them provides a proof of the Main Theorem.

In \S\ref{S:connections} we describe the most important metric connections in the tangent bundle $TM$ and show that most of the terms in the Allendoerfer-Weil formula vanish for the Levi-Civita connection.  In \S \ref{S:proof}, we describe the process of averaging over metric connections, and use this averaging procedure to average over the action of the symmetry group $G$, thus obtaining an Allendoerfer-Weil formula for which all terms vanish except those in dimensions zero and either $n$ when $n$ is even or $n-1$ when $n$ is odd.  This finally enables us to finish the proof of the Hazzidakis formulae and the Main Theorem.

As a graduate student, the author admired the beautiful extensions of the Gauss-Bonnet formula for surfaces that had been developed by Allendoerfer, Fenchel, Weil and Chern.  Our conjectured generalization of Hilbert's theorem resulted from a search for a nice geometric application of these formulae to submanifold theory.  More generally, we were curious as to whether the representation of characteristic classes via differential forms might yield nontrivial global theorems relating curvature to topology for $n$-dimensional submanifolds of low codimension in Euclidean space.  The proof of the Main Theorem does vindicate the author's suspicion based on early calculations that the conjecture could be solved via the Gauss-Bonnet formula.  The author is indebted to Shoshichi Kobayashi (a student of Allendoerfer), Hung-hsi Wu and Shiing-shen Chern for helping him embark on a longer than expected journey.
 
\section{The principal adapted moving frame}
\label{S:adaptedframe}

We begin by reviewing the local theory of submanifolds of constant negative curvature using Cartan's method of moving frames, as presented in \cite{Mo}, with a few minor changes in the arguments.

We can study an $n$-dimensional oriented Riemannian manifold $M$ by means of an oriented moving frame, an ordered $n$-tuple of orthogonal unit-length vector fields $(e_1, \ldots , e_n)$ defined on an open subset $U$ of $M$ with dual one-forms $(\theta _1, \ldots , \theta_n)$ defined by
$$\theta _i(e_j) = \delta _{ij} = \begin{cases} 1, & \hbox{if $i=j$}, \cr 0, & \hbox{if $i\neq j$,}\end{cases}$$
which satisfy the additional requirement that $(e_1, \ldots , e_n)$ be positively oriented, that is, that
$$\theta _1 \wedge \theta _2 \wedge \cdots \wedge \theta _n$$
agrees with the positively oriented volume form on $M$.  The {\em Levi-Civita connection\/} is defined via the corresponding connection one-forms $\omega _{ij}$, for $1 \leq i,j \leq n$, which satisfy the structure equations
$$d\theta _i = \sum _{j=1}^n \omega _{ij} \wedge \theta _j, \quad \omega _{ij} + \omega _{ji} = 0.$$
These connection forms determine the covariant differential $\nabla V$ of a vector field $V$ tangent to $M$ by the formula
\begin{equation} V = \sum_{i=1}^nv_i e_i \quad \Rightarrow \quad \nabla V = \sum_{i=1}^n \left(dv_i + \sum_{j=1}^n \omega _{ij} v_j \right) e_i. \label{E:covdiff} \end{equation}
The curvature of the Levi-Civita connection is expressed by the family of two-forms
$$\Omega _{ij} = d \omega _{ij} + \sum _{j=1}^n \omega _{ik} \wedge \omega _{kj}, \quad \hbox{which satisfy} \quad \Omega _{ij} + \Omega _{ji} = 0.$$
The fact that $M$ has constant curvature $k = -1$ is then expressed by the equation
$$\Omega _{ij} = k \ \theta _i \wedge \theta _j = - \ \theta _i \wedge \theta _j.$$

Throughout this article, ${\mathbb H}^n$ will denote $n$-dimensional hyperbolic space, the unique complete simply connected $n$-dimensional Riemannian manifold of constant curvature $-1$, and ${\mathbb E}^N$ will denote $N$-dimensional Euclidean space.  We suppose that $M^n$ is an oriented simply connected Riemannian manifold of constant curvature $-1$, not necessarily complete, and that we are given an isometric immersion $f: M^n \longrightarrow {\mathbb E}^{2n-1}$.  If $h: {\mathbb E}^{2n-1} \longrightarrow {\mathbb H}^{2n}$ is the standard isometric imbedding from ${\mathbb E}^{2n-1}$ onto a horosphere, then the composition $g = h \circ f$ can be regarded as a \lq\lq developable submanifold\rq\rq \ of ${\mathbb H}^{2n}$.  We will use the following conventions on the ranges of indices:
$$1 \leq i,j,k \leq n; \qquad n + 1 \leq \lambda ,\mu ,\nu \leq 2n; \qquad 1 \leq I, J, K \leq 2n.$$

We extend the moving frame above to an adapted moving frame over a contractible open subset $U$ of $M^n$, a collection $(e_1 , \ldots , e_{2n})$ of sections of the pullback $g^*T ({\mathbb H}^{2n})|U$ of the tangent bundle to ${\mathbb H}^{2n}$ such that:
\begin{enumerate}
\item $e_1 , \ldots , e_{2n}$ are orthonormal with respect to the Riemannian metric $\langle \cdot, \cdot \rangle$ on ${\mathbb H}^{2n}$,
\item $(e_1, \ldots , e_n)$ restrict to a positively oriented moving frame for the oriented Riemannian manifold $U \subseteq M$, and hence the remaining vectors are normal to $g(U)$.
\item $e_{2n}$ is perpendicular to the horosphere $h({\mathbb E}^{2n-1})$.
\end{enumerate}
We extend the coframe $(\theta _1, \ldots , \theta _n)$, defined by $\theta _i (e_j) = \delta _{ij}$, to $(\theta _{1}, \ldots , \theta _{2n})$ by setting $\theta _\lambda = 0$ on $U$.  The Levi-Civita connection on ${\mathbb H}^{2n}$ induces a connection $\nabla $ on $g^*T({\mathbb H}^{2n})|U$, which has corresponding connection one-forms $\omega _{IJ}$ defined by
$$\omega _{IJ}= \langle e_I, \nabla e_J \rangle, \quad \omega _{IJ} + \omega _{JI} = 0.$$
These differential forms satisfy the Cartan structure equations
\begin{equation} d\theta _I = - \sum \omega _{IJ} \wedge \theta _J, \label{eq:hilb21} \end{equation}
\begin{equation} d\omega _{IJ} = - \sum \omega _{IK} \wedge \omega _{KJ} -\theta _I \wedge \theta _J, \label{eq:hilb22} \end{equation}
the second of these expressing the fact that ${\mathbb H}^{2n}$ has constant curvature $-1$.

It follows from $\theta _\lambda = 0$ and equation (\ref{eq:hilb21}) that
$$\sum \omega _{\lambda i} \wedge \theta _i = 0 \qquad \hbox{and hence} \qquad \omega _{\lambda i} = \sum h_{\lambda ij} \theta _j,$$
where the $h_{\lambda ij}$\rq s are smooth functions on $U$ such that $h_{\lambda ij} = h_{\lambda ji}$.  The symmetric bilinear forms
$$\Phi _\lambda = \sum \omega _{\lambda i} \otimes \theta _i = \sum h_{\lambda ij}\theta _i \otimes \theta _j$$
are called the {\em second fundamental forms\/} of $g$.  The assumption that the sectional curvatures of $M^n$ are $-1$ implies that
$$d\omega _{ij} = - \sum \omega _{ik} \wedge \omega _{kj} - \theta _i \wedge \theta _j,$$
and comparison with (\ref{eq:hilb22}) yields
$$\sum \omega _{\lambda i} \wedge \omega _{\lambda j} = 0,$$
or equivalently,
\begin{equation} \sum [\Phi _\lambda (x,z) \Phi _\lambda (y,w) - \Phi _\lambda (x,w) \Phi _\lambda (y,z)] = 0, \label{eq:hilb23} \end{equation}
for $x,y,z,w \in T_pM$.  Following the terminology of \'Elie Cartan \cite{Ca}, we say that (\ref{eq:hilb23}) expresses the fact that the second fundamental forms of $g: M^n \longrightarrow {\mathbb H}^{2n}$ are {\em exteriorly orthogonal\/}.  The condition that $e_{2n}$ is perpendicular to the horosphere implies that
\begin{equation}\omega _{2n,i} = \theta _i, \qquad \hbox{and hence} \qquad \Phi _{2n} = \sum \theta _i \otimes \theta _i.\label{eq:hilb24}\end{equation}

\vskip .1 in
\noindent
{\bf Lemma~2.1.} {\sl There exist smooth one-forms $\phi_1, \ldots \phi _n$ on $U$ and an $n\times n$ matrix of 
smooth real-valued functions $A = (a_{ij})$ such that $A$ is orthogonal and
\begin{equation} \Phi _{n+i} = \sum a_{ij} \phi _j \otimes \phi _j. \label{eq:hilb25}\end{equation}
}

\noindent
This Lemma follows from an algebraic result of Cartan \cite{Ca} which we described and proved in \cite{Mo}.  Our earlier argument can be made a little simpler by using the theory of flat bilinear forms developed in \cite{Mo77a} and \cite{Mo77}.

If $(V, \langle , \rangle)$ and $(W, \langle , \rangle)$ are real inner product spaces of dimension $n$ and
$$\beta : V \times V \longrightarrow W$$
is a symmetric bilinear map such that
\begin{equation} \langle \beta (x,x),  \beta (y,y) \rangle - \langle \beta(x,y),  \beta(x,y) \rangle = 0, \quad \hbox{for} \quad x, y \in V, \label{E:flatid}\end{equation}
we say that $\beta $ is {\em flat\/} with respect to the inner product $\langle \cdot , \cdot \rangle$ on $V$.  Let
$$N(\beta ) = \{ x \in V : \hbox{$\beta(x,y) = 0$ for all $y \in V$}.$$

The example we have in mind is
\begin{equation} V = T_pM, \quad W = N_pM, \quad \beta = \sum _{\lambda = n+1}^{2n} \Phi _\lambda e_\lambda, \label{E:exa}\end{equation}
where $(e_{n+1}, \ldots, e_{2n})$ is an orthonormal basis for the normal space $N_pM$.  For this example, $N(\beta ) = 0$ by (\ref{eq:hilb24}).  If $x \in V$ we define $\beta (x) : V \rightarrow W$ by $\beta (x)(y) = \beta (x,y)$ and let
$$q = \max \{ \dim \beta (x)(V) : x \in V \}.$$
We can then reformulate Cartan's result as follows:

\vskip .1 in
\noindent
{\bf Lemma~2.2.} {\sl If $(V, \langle , \rangle)$ and $(W, \langle , \rangle)$ are $n$-dimensional real inner product spaces and
$$\beta : V \times V \longrightarrow W$$
is a flat bilinear form with $N(\beta ) = 0$, then there exists an orthonormal basis $(v_1, \ldots, v_n)$ of $V$ such that
\begin{equation} \dim \beta (v_i)(V) = 1, \quad \hbox{for $1 \leq i \leq n$.} \label{E:induct}\end{equation}}

\noindent
Proof:  We first adapt the proof of the Lemma in \S 2 of \cite{Mo77a}.  Let
$$q = \max \{ \dim \beta (x)(V) : x \in V \}.$$
We say that $x \in V$ is a {\em regular element \/} if $\dim \beta (x)(V) = q$.  The regular elements form an open dense subset of $V$.

We now let
$$N(\beta,x) = \{ n\in V : \beta (x)(n) = 0 \},$$
and claim that if $x$ is a regular element and $n \in N(\beta ,x)$,
\begin{equation} \langle \beta (n,y), \beta (x,V) \rangle = 0 \quad \hbox{and} \quad \beta (n,y) \in \beta (x,V) \quad \hbox{for all $y \in V$.}\label{E:claim}\end{equation}
The first of these assertions is a consequence of (\ref{E:flatid}):
$$z \in V \quad \Rightarrow \quad \langle \beta (n,y), \beta (x,z) \rangle = \langle \beta (n,x), \beta (y,z) \rangle = 0.$$
To prove the second, we note that $\beta (x + ty)(V)$ is a continuously varying $q$-dimensional subspace of $W$ for $t$ sufficiently small since $x$ is regular.  Then
$$\beta (x+ty,n) = t \beta (y,n) \Rightarrow \beta (y,n) \in \beta (x+ty)(V) \quad \hbox{for all $t$},$$
and when we take the limit as $t \rightarrow 0$ we find that $\beta (n,y) \in \beta (x,V)$.

Since $N(\beta ) = 0$ and the inner product on $W$ is positive definite, (\ref{E:claim}) implies that when $x$ is regular, $N(\beta, x) =0$ and $\beta (x) : V \rightarrow W$ is an isomorphism.

We next adapt the argument from Theorem 2 of \S 5 from \cite{Mo77}.  We choose a regular element $x \in V$ with $\beta(x)$ being an isomorphism and for $y \in V$, define
$$B(y): W \rightarrow W \quad \hbox{by} \quad B(y) = \beta (y) \circ \beta (x)^{-1}.$$
Suppose that if $v_1, v_2 \in V$ and
$$w_1 = \beta(x)(v_1), \quad w_2 = \beta (x) v_2.$$
We then use (\ref{E:flatid}) to show that
$$\langle B(y)w_1, w_2 \rangle = \langle \beta (y) v_1, \beta(x) v_2 \rangle = \langle \beta (y) v_2, \beta(x) v_1 \rangle = \langle B(y)w_2, w_1 \rangle, $$
so $B(y)$ is symmetric with respect to the inner product on $W$.  Again we use (\ref{E:flatid}), this time to show also that if $y_1,y_2 \in V$, then
\begin{multline*} \langle B(y_1)B(y_2) w_1, w_2 \rangle = \langle B(y_2) w_1, B(y_1)w_2 = \langle \beta (y_2)(v_1), \beta (y_1)(v_2) \rangle \\
= \langle \beta (y_1)(v_1), \beta (y_2)(v_2) \rangle = \langle B(y_1) w_1, B(y_2)w_2 \rangle = \langle B(y_2)B(y_1) w_1, w_2 \rangle, \end{multline*}
so $B(_1)$ and $B(_2)$ commute.  Thus the $\{ B(y) : y \in V\}$ form a commuting family of linear endomorphisms of the inner product space $(W, \langle \cdot , \cdot \rangle )$ which are symmetric with respect to the inner product, and can therefore be simultaneously diagonalized by an orthonormal basis for $(W, \langle \cdot , \cdot \rangle )$.  Note that this $n$-dimensional family is spanned by matrices of rank one.

Finally, we establish (\ref{E:induct}) by induction on the dimension $n$ of $V$ and $W$:  The assertion is clearly true when $n = 1$.  To establish the inductive step when $\dim V=\dim W = n$, we choose an element $v_n \in V$ such that $\dim \beta (v_n)(V) = 1$ and let
$$V_n = N(\beta ,v_n), \quad W_n = \hbox{orthogonal complement of $\beta (v_n)$.}$$
Then the fact that $\beta $ is flat implies that it induces a flat bilinear map
$$\beta_n : V_n \times V_n \longrightarrow W_n,$$
and by inductive hypothesis $V_n$ possesses an orthonormal basis $(v_1, \ldots, v_{n-1})$ with
$$\dim \beta_0(v_i) = 1, \quad \hbox{for $1 \leq i \leq n-1$.}$$
We can then add $v_n$ to achieve an orthonormal basis for $V$ satisfying (\ref{E:induct}).QED 

\vskip .1 in
\noindent
When we apply Lemma~ 2.2 to our example (\ref{E:exa}), we immediately obtain the conclusion of Lemma~2.1.  We can rewrite (\ref{eq:hilb25}) as
\begin{equation}\sum a_{ji} \Phi _{n+i} = \phi_j \otimes \phi_j.\label{eq:hilb26}\end{equation}
Thus if we relax the condition that $e_{2n}$ be perpendicular to the horosphere, we can construct a new moving orthonormal frame $(e_{n+1},\ldots ,e_{2n})$ for the normal bundle such that
$$\Phi _{n+i} = \sum \omega _{n+i,j} \otimes \theta _j = \phi _i \otimes \phi _i = \frac{1}{x_i} \theta _i \otimes \theta _i.$$
The corresponding new adapted orthonormal moving frame $(e_1,\ldots ,e_{2n})$ will be called the {\em principal adapted frame\/}.  Although we have only constructed this frame locally, uniqueness up to changes of sign and permutations implies that this principal adapted frame can be defined globally, if $M$ is simply connected.

Equation (\ref{eq:hilb26}) implies that
\begin{equation}\omega _{n+i,j} = \frac{\delta _{ij}}{ x_i} \theta _i.\label{eq:hilb27}\end{equation}
We differentiate this equation when $i=j$ to obtain
$$d\omega _{n+i,i} = - \frac{dx_i}{x_i^2} \wedge \theta _i - \frac{1}{x_i} \sum \omega _{ij} \wedge \theta _j.$$
On the other hand, it follows from $\theta _\lambda = 0$, (\ref{eq:hilb22}) and (\ref{eq:hilb27}) that
$$d\omega _{n+i,i} = - \sum \omega _{n+i,j} \wedge \omega _{ji} - \sum \omega _{n+i,n+j} \wedge \omega _{n+j,i} = 0.$$
Comparing these two equations, we conclude that
$$\frac{1}{x_i^2}e_j(x_i) = \frac{1}{x_i} \omega _{ij}(e_i),$$
where $e_j(x_i)$ denotes the \lq\lq directional derivative\rq\rq \ of $x_i$ in the direction of $e_j$.  It follows immediately that
\begin{equation}\omega _{ij} = \frac{1}{x_i}e_j(x_i)\theta _i - \frac{1}{x_j}e_i(x_j)\theta _j. \label{eq:hilb28}\end{equation}

Differentiation of equation (\ref{eq:hilb27}) in the case where $i \neq j$ yields
$$0 = d\omega _{n+i,j} = \omega _{n+i,i}\wedge \omega _{ij} - \omega _{n+i,n+j} \wedge \omega _{n+j,j} \qquad \qquad $$
$$\qquad \qquad = - \frac{1}{x_i}\theta _i \wedge \left[ \frac{1}{x_i}e_j(x_i)\theta _i - \frac{1}{ x_j}e_i(x_j)\theta _j \right] - \omega _{n+i,n+j} \wedge \frac{1}{x_j}\theta_j$$
$$= \frac{1}{x_ix_j} e_i(x_j) \theta _i \wedge \theta _j - \frac{1}{x_j}\omega _{n+i,n+j}\wedge \theta _j.$$
We conclude that $\omega _{n+i,n+j}(e_k) = 0$ unless $k = i$ or $j$ and
$$\omega _{n+i,n+j}(e_i) = \frac{1}{x_i} e_i(x_j), \quad \omega _{n+i,n+j}(e_j) = - \frac{1}{x_j} e_j(x_i).$$
We thus derive the formula,
\begin{equation}\omega _{n+i,n+j} = \frac{1}{ x_i} e_i(x_j) \theta _i - \frac{1}{x_j} e_j(x_i) \theta _j. \label{eq:hilb29}\end{equation}

A calculation based on (\ref{eq:hilb27}) quickly shows that
$$d\left( \frac{1}{x_i} \theta _i \right) = 0,$$
and hence our standing assumption that $M$ is simply connected implies that once a base point $p \in M$ is chosen, there exist unique smooth functions $y_i:M \rightarrow {\mathbb R}$ such that
\begin{equation}dy_i = \frac{1}{x_i} \theta _i, \qquad y_i(p) = 0.\label{eq:hilb210}\end{equation}

\vskip .1 in
\noindent
{\bf Lemma~2.3.} {\sl If $M$ is complete as well as simply connected, $\psi = (y_1, \ldots , y_n)$ is a global coordinate system on $M$.}
\vskip .1 in

\vskip .1 in
\noindent
Let us recall the key ideas of the argument we gave in \cite{Mo}.  The vector fields dual to $dy_1, \ldots , dy_n$ are
$$Y_1 = x_1 e_1, \quad \ldots , \quad Y_n = x_n e_n.$$
Since these vector fields have length $< 1$, they are complete, and possess global one-parameter groups of diffeomorphisms.  Define a map $\phi : {\mathbb R}^n \rightarrow M$ by
$$\phi (t_1,\ldots ,t_n) = \phi_1(t_1,\phi_2(t_2, \ldots ,\phi_n(t_n,p) \ldots )),$$
where $\phi _i$ is the one-parameter group of $Y_i$.  It follows from (\ref{eq:hilb210}) that $[Y_i,Y_j] = 0$, for all $i,j$, so that the one-parameter groups commute, and hence
$$\phi _*(\partial /\partial t_i ) = Y_i.$$
To finish the proof of the Lemma, one needs only show that
$$\psi \circ \phi = \hbox{id}, \qquad \phi \circ \psi = \hbox{id},$$
which is easily done by showing that the sets on which these equations hold are both open and closed.  QED

\vskip .1 in
\noindent
The functions $y_1$, \ldots $y_n$ we have constructed define what is called a {\em principal coordinate system\/} centered at the point $p \in M$.  If $(y_1, \ldots , y_n)$ and $(\hat y_1, \ldots , \hat y_n)$ are two such principal coordinate systems centered at different points of $M$ with respect to the same principal adapted frame, then
$$y_i - \hat y_i = c_i, \quad \hbox{for $1 \leq i \leq n$,}$$
where the $c_i$'s are suitable constants.

The unit-length vector fields
\begin{equation} \pm x_1e_1 \pm \cdots \pm x_ne_n \label{E:asymptoticvectorfields} \end{equation}
are called {\em asymptotic vector fields\/} on $M$.  There are $2^n$ such asymptotic vector fields corresponding to the $2^n$ possible choices of signs.  The functions 
$$\pm y_1 \pm \cdots \pm y_n$$
are called {\em asymptotic coordinates\/}.  A choice of $n$ linearly independent asymptotic coordinates $(z_1, \ldots , z_n)$ defines a corresponding {\em asymptotic coordinate system\/} on $M$ whose coordinate vectors are unit-length asymptotic vector fields on $M$, which we denote by $Z_1, \ldots , Z_n$.  The construction of such principal and asymptotic coordinate systems completes our extension of the first step in the classical proof of Hilbert's theorem for surfaces in ${\mathbb E}^3$.

Choice of a distinguished asymptotic vector field
\begin{equation} X = x_1 e_1 + \cdots + x_n e_n \label{E:asymvect} \end{equation}
allows us to construct a vector bundle isomorphism
$$F_X : TM \rightarrow N(h \circ f) = (\hbox{normal bundle of $h \circ f$}), \quad F_X(v) = \sum e_\lambda \Phi _\lambda (X,v),$$
which we use to pull the metric connection in the normal bundle back to a metric connection $\widehat \nabla $ on $TM$ (which does not have vanishing torsion).  Note that $F_X$ takes our distinguished asymptotic vector field to a normal vector field
$$\widehat X = x_1e_{n+1} + \cdots + x_ne_{2n},$$
which is orthogonal to the horosphere ${\mathbb E}^{2n-1} \subseteq {\mathbb H}^{2n}$ and is therefore parallel with respect to the connection in the normal bundle.  This implies that $X$ is parallel with respect to the pullback connection $\widehat \nabla $.  One can check that since the second fundamental forms $\Phi _\lambda $ can be diagonalized simultaneously, this connection $\widehat \nabla $ is flat.  We will sometimes use the notation
\begin{equation} \phi _{ij} = \omega _{n+i,n+j}, \qquad \hbox{so that} \qquad \widehat \nabla e_j = \sum e_i \phi _{ij}. \label{E:flatconnection} \end{equation}

\vskip .1 in
\noindent
{\bf Dual case: submanifolds of constant positive curvature.}  It is a curious fact that the theory we have just presented has an analog for n-dimensional submanifolds of constant positive curvature in ${\mathbb E}^{2n-1}$, which is described in \cite{Mo77}.  If $f : U \rightarrow {\mathbb E}^{2n-1}$ is an isometric immersion from an open subset $U$ of the unit curvature sphere ${\mathbb S}^{n}$, a point $p \in U$ is said to be a {\em weak umbilic\/} if a neighborhood of $p$ has the form expected for a composition of $U \rightarrow {\mathbb E}^{n+1}$ with a developable immersion from a neighborhood of $f(p)$ into ${\mathbb E}^{2n-1}$.  Otherwise, it is said to be a {\em nonumbilic\/}.  There is a rich supply of nonumbilic submanifolds of constant positive curvature which correspond to solutions to a generalization of a sinh-Gordon equation which can be regarded as dual to the sine-Gordon equation.  But it is proven in \cite{Mo77} that nonumbilics cannot occur for global isometric immersions of ${\mathbb S}^{n}$ in ${\mathbb E}^{2n-1}$.  The author conjectures that this result implies that any smooth isometric immersion $f: S^n \rightarrow {\mathbb E}^N$, where $S^n$ is the round sphere of constant curvature one and $N \leq 2n-1$, can be extended to a smooth isometric immersion of the flat $(n+1)$-disk $D^{n+1}$ which has boundary $S^n$, a fact proven in \cite{Mo96} under the additional condition that $N \leq n+2$.

\section{Asymptotic polyhedra}
\label{S:asymptotic}

In the even-dimensional case, the proof of the Main Theorem will be based upon estimating the volume of a family of {\em asymptotic polyhedra\/} defined in terms of a principal coordinate system $(y_1, \ldots y_n)$ on a hyperbolic space ${\mathbb H}^n$ which is isometrically immersed in ${\mathbb R}^{2n-1}$ by 
\begin{equation} P_a = \{ p \in {\mathbb H}^n : |\pm y_1(p)\pm \cdots \pm y_n(p) | \leq a, \hbox{ for all choices of signs}\}. \label{E:Riemannianpolyhedron}\end{equation}
for $a \in (0, \infty)$.  In the odd-dimensional case, we estimate instead the volume of the boundaries $\partial P_a$ of these asymptotic polyhedra.

We regard $P_a$ as centered at the point $p_0 \in {\mathbb H}^n$ with coordinates
$$y_1 = y_2 = \cdots = y_n = 0.$$
For a point $p$ of ${\mathbb H}^n$ to lie within $P_a$, the value of each asymptotic coordinate
$$z = \pm y_1 \pm y_2 \pm \cdot \pm y_n$$
at $p$ must have absolute value $\leq a$, for all $2^n$ possible choices of sign.  The asymptotic polyhedron $P_a$ specializes to an asymptotic rectangle (\ref{E:asymptoticrectangle}) of the form considered in the Hazzidakis formula for surfaces in the case where the sides have equal coordinate length.

Note that the boundary $\partial P_a$ has a simplicial decomposition with $2n$ zero-simplices located along the coordinate axes and $2^n$ $(n-1)$-simplices.  Each of these $(n-1)$-simplices is associated to a unique unit-length asymptotic vector field which is transverse to that simplex and points outward from $p_0$.

This polyhedron $P_a$ has a large symmetry group, the group $G$ generated by reflections in the principal coordinate hyperplanes
$$y_1 = 0, \quad y_2 = 0, \quad \ldots , \quad y_n = 0,$$
each such reflection being defined by the explicit formula
\begin{equation} g_i(y_j) = (1 - 2 \delta _{ij}) y_j, \quad \hbox{for $1 \leq j \leq n$.} \label{E:action} \end{equation}
Note, however, that this action is not usually an isometric action on $P_a$ in terms of the Riemannian metric on ${\mathbb H}^n$.  On the other hand, there is a related action of $G$ as isometries on each fiber of the tangent bundle $T({\mathbb H}^n)$.  In terms of the principal frame $(e_1, e_2, \ldots , e_n)$, the operation of the generator $g_i \in G$ on $T_p({\mathbb H}^n)$ is given by
$$g_i(e_i(p)) = - e_i(p) \quad \hbox{and} \quad g_i(e_j(p)) = e_j(p) \quad \hbox{for $j \neq i$.}$$

The group $G$ is abelian of order $2^n$, and each element $g \in G$ corresponds to a map
$$g : \{ 1 , \ldots ,n \} \rightarrow \{ -1,1 \} \qquad  \hbox{such that} \qquad g(y_i) = g(i) y_i.$$
The group $G$ acts simply transitively on the unit-length asymptotic vector fields, the element $g \in G$ sending $X$ defined by (\ref{E:asymvect}) to
$$X_g = g(1) x_1 e_1 + g(2) x_2 e_2 + \cdots + g(n) x_n e_n.$$
The group $G$ also acts simply transitively on the $(n-1)$-simplices in the simplicial decomposition of $\partial P_a$, which are in one-to-one correspondence with the asymptotic vector fields.

We let
\begin{equation} \epsilon (g) = g(1) \cdots g(n) = \left\{\begin{array}{ll} 1 & \mbox{if $g$ is orientation-preserving,} \\ -1 & \mbox{if $g$ is orientation-reversing,} \end{array} \right. \label{E:epsilong} \end{equation}
and define the subgroup
$$G_0 = \{g \in G : \epsilon (g) = 1 \},$$
of index two in $G$, which consists of the orientation-preserving elements.  Note that elements of $G_0$ lie in the rotation group $SO(n)$ while elements of $G$ lie in its double cover $\hbox{Spin}(n)$.  Since elements of $G_0$ preserve the orientation of $T({\mathbb H}^n)$, they will preserve the Euler class of $T({\mathbb H}^n)$ that will be defined in the next section.

Finally, $G$ acts on the space of metric connections on $T({\mathbb H}^n)$, which are not required to have vanishing torsion, this action
$$\nabla \mapsto \nabla ^g \quad \hbox{being defined by} \quad \nabla ^g = g \circ \nabla \circ g, \quad \hbox{for $g \in G$.}$$
(If $\nabla $ is the Levi-Civita connection, $\nabla ^g$ might not have vanishing torsion because $g$ is not usually the differential of an isometric action on $U \subseteq {\mathbb H}^n$.)

Our formula for the volume of $P_a$ or $\partial P_a$ will follow from the generalized Gauss-Bonnet formula for $P_a$, a manifold with boundary that is not smooth, but has corners.  This formula is an adaptation of the formula found by Allendoerfer and Weil.  Eventually we will take the signed average of the Allendoerfer-Weil formulae for an asymptotic polyhedron $P_a$ over the collection of connections $\{ \nabla ^g : g \in G \}$.  We do this by taking the signed average of each term in the Allendoerfer-Weil formula, the sign being adjusted by multiplying each summand by the expression $\epsilon (g)$ of (\ref{E:epsilong}).  We need to take the signed average because we fix the orientation on the base $P_a$ while the action of $G$ can reverse orientation on the fiber of $TM$.

\section{The Gauss-Bonnet Theorem}
\label{S:gaussbonnet}

If $A = (a_{ij})$ is a skew-symmetric matrix with $2m$ rows and $2m$ columns, the Pfaffian of $A$ is
\begin{equation} \hbox{Pf}(A) = \frac{1}{2^mm!}\sum _{\sigma \in S_{2m}} \hbox{sgn}(\sigma ) a_{\sigma (1),\sigma (2)} \cdots a_{\sigma (2m-1)\sigma (2m)}. \label{E:pfaff}\end{equation}
Here $S_{2m}$ is the group of permutations on $2m$ letters and $\hbox{sgn}(\sigma )$ is the sign of the permutation $\sigma \in S_{2m}$.  The factor $(1/2^m)m!$ is included in the definition so that
$$\hbox{Pf}\begin{pmatrix} 0 & a_1 & 0 & 0 & \cdots & 0 & 0 \cr
-a_1 & 0 & 0 & 0 & \cdots & 0 & 0 \cr
0 & 0 & 0 & a_2 & \cdots & 0 & 0 \cr
0 & 0 & -a_2 & 0 & \cdots & 0 & 0 \cr
\cdot & \cdot & \cdot & \cdot & \cdots & \cdot & \cdot \cr
0 & 0 & 0 & 0 & \cdots & 0 & a_m \cr
0 & 0 & 0 & 0 & \cdots & -a_m & 0 \cr \end{pmatrix} = a_1 a_2 \cdots a_m.$$
One can verify that if $B \in SO(2m)$, then
\begin{equation} \hbox{Pf}(BAB^{-1}) = \hbox{Pf}(A) \label{E:pfaffian} \end{equation}
Moreover, $\hbox{Pf}(A)^2 = \det (A)$.  Thus the Pfaffian serves as a square root of the determinant for skew-symmetric matrices.

If $M$ is a $(2m)$-dimensional Riemannian manifold and $\Omega = (\Omega _{ij})$ is the matrix of curvature two-forms corresponding to a positively oriented moving frame for the tangent bundle of $M$, we can construct a form $\hbox{Pf}(\Omega )$ of degree $2m$ on $M$.  It follows from (\ref{E:pfaffian}) that $\hbox{Pf}(\Omega )$ does not depend on the choice of positively oriented moving orthonormal frame, and is therefore a globally defined $(2m)$-form on $M$.  One easily checks that if the metric has constant curvature $-1$ as in \S \ref{S:adaptedframe},
\begin{multline} (-1)^m\hbox{Pf}(\Omega ) = c_m \theta _1 \wedge \theta _2 \wedge \cdots \wedge \theta _{2m}  = c_m (\hbox{volume form on $M$}), \\ \hbox{where} \quad c_m = \frac{(2m)!}{2^m m!} = (2m-1) \cdot (2m-3) \cdots 3 \cdot 1. \label{E:pf=vol}\end{multline}

\vskip .1 in
\noindent
{\bf Gauss-Bonnet Theorem.} {\sl If $M$ is an oriented compact $(2m)$-dimensional Riemannian manifold, then
$$\frac{1}{(2\pi)^m} \int _M \hbox{Pf}(\Omega ) = \chi (M) = (\hbox{Euler characteristic of $M$}),$$
where $\Omega $ is the curvature of the Levi-Civita connection in the oriented tangent bundle $TM$.  If $M$ is an oriented compact odd-dimensional manifold, $\chi (M) = 0$.}

\vskip .1 in
\noindent
We recall some well-known history (see \cite{KN}, pages 357-361).  The Gauss-Bonnet Theorem was first proven by Allendoerfer \cite{A} and Fenchel \cite{F} under the assumption that $M$ admits an isometric imbedding into some Euclidean space ${\mathbb E}^N$.  This would have proven the theorem for general Riemannian manifolds if Nash's isometric imbedding theorem had been known.   Allendoerfer and Weil \cite{AW} next gave a proof for a compact \lq\lq Riemannian polyhedron" which is isometrically imbedded in ${\mathbb E}^N$.  Since it was known that real analytic Riemannian manifolds admit \lq\lq local" isometric imbeddings into Euclidean space by a theorem of Janet and Cartan \cite{Ca1}, one could then piece together the result for Riemannian polyhedra to give a proof for all real analytic Riemannian manifolds.  Then Chern \cite{Ch} gave a proof of the Gauss-Bonnet Theorem for smooth compact Riemannian manifolds which avoided the need of an isometric imbedding and illustrated the power of Cartan's method of moving frames.

Chern's argument led to the theory of characteristic classes as described in Appendix~C of Milnor and Stasheff \cite{MS}, specialized to the Euler class of an oriented geometric vector bundle over $M$.  An {\em oriented geometric vector bundle\/} over a smooth manifold $M$ is a triple $(E,\langle , \rangle, \nabla)$, with $E$ being an oriented real vector bundle over $M$, $\langle , \rangle $ a fiber metric on $E$ and $\nabla $ a connection in $E$ which preserves the fiber metric.  If $E$ has rank $n$, the structure group of the oriented geometric vector bundle is $SO(n)$ and the Pfaffian is an \lq\lq invariant polynomial in the curvature $\Omega $" associated to $SO(n)$ when $n$ is even.  The Gauss-Bonnet Theorem extends to say that if $(E,\langle , \rangle, \nabla)$ is an oriented geometric vector bundle of rank $n$ over the compact $n$-dimensional manifold $M$, and
\begin{equation} e(\nabla ) = \begin{cases} \frac{1}{(2\pi)^m}\hbox{Pf}(\Omega ), & \hbox{when $n = 2m$ is even,} \cr
0,& \hbox{when $n$ is odd.} \end{cases} \label{E:eulerform}\end{equation}
where $\Omega $ denotes the curvature of $\nabla $, then $\int _M e(\nabla )$ is a topological invariant of $E$, hence independent of metric connection, and this invariant equals the Euler characteristic $\chi(M)$ when $E = TM$.  For example, we can apply the Gauss-Bonnet formula to the connections $\nabla ^g$ obtained from the Levi-Civita connection $\nabla $ by the action of $G_0$ that were described in the previous section.

Later Matthai and Quillen (\cite{Q} and \cite{MQ}) gave a proof using the theory of superconnections that clarifies the role of the {\em Thom form\/} of the geometric vector bundle, a smooth $n$-form on $E$ that is rapidly decreasing in the fiber direction.  We can adapt their approach as follows.

The simplest case is that of the trivial bundle $M \times {\mathbb R}^n$, with standard fiber metric and trivial connection.  In this case, we take positively oriented Euclidean coordinates $(t_1, t_2, \ldots , t_n)$ on the fiber ${\mathbb R}^n$, and define the {\em Thom form\/} on the trivial bundle $M \times {\mathbb R}^n$ by
\begin{equation} \tau = \left( \frac{1}{\sqrt \pi }\right) ^n e^{-(t_1^2 + \cdots + t_n^2)} dt_1 \wedge \cdots \wedge dt_n. \label{E:flatthom}\end{equation}
This form is closed and invariant under the action of the special orthogonal group on the fiber, and since
\begin{equation}\int _{-\infty}^\infty e^{-t^2}dt = \sqrt \pi, \label{E:gaussian}\end{equation}
it has the property that its integral over the fiber ${\mathbb R}^n$ is one.  We can also write the Thom form in polar coordinates as
$$\tau = \left( \frac{1}{\sqrt \pi }\right) ^n e^{-r^2} r^{n-1}dr \wedge \Theta,$$
where $\Theta $ is the volume form on the unit $(n-1)$-sphere.  If we set $u = r^2$,
$$\int _0^\infty  e^{-r^2} r^{n-1} dr = \frac{1}{2} \int _0^\infty e^{-u} u^{(n/2) - 1} du = \frac{1}{2} \Gamma\left(\frac{n}{2}\right),$$
by definition of the gamma function.  Integration over the fiber in polar coordinates therefore yields a well-known formula for the volume of the $(n-1)$-sphere,
\begin{equation} \hbox{Volume of $S^{n-1}$} = \frac{2(\sqrt \pi)^n}{\Gamma (n/2)},\label{E:volSn-1}\end{equation}
which specializes when $n$ is even, say $n = 2m$, to
\begin{equation} \hbox{Volume of $S^{2m-1}$} = \frac{2\pi^m}{(m-1)!}.\label{E::volS2m-1}\end{equation}

To construct the Thom form for a general geometric vector bundle $(E,\langle , \rangle, \nabla)$ we add a curvature term to (\ref{E:flatthom}).  To construct this systematicaly, we consider the algebra of differential forms with values in the exterior algebra of $E$:  Let
$${\mathcal A} ^{p,q} ={\mathcal A} ^p(M,\Lambda ^q(E)) = \{ \hbox{smooth $p$-forms on $M$ with values in $\Lambda ^q(E)$} \},$$
and give the direct sum ${\mathcal A} ^{*,*} = \sum _{p,q} {\mathcal A} ^{p,q}$ the skew-commutative wedge product
$$(\phi \alpha ) \wedge (\psi \beta ) = (-1)^{\deg \psi \deg \alpha }(\phi \wedge \psi )(\alpha \wedge \beta),$$
for $\phi ,\psi \in {\mathcal A} ^{*,0}$, $\alpha ,\beta \in{\mathcal A} ^{0,*}$.  This makes ${\mathcal A} ^{*,*}$ into a skew-commutative associative algebra with unit over the ring of functions.  The connection $\nabla $ extends to a family of linear maps
$$\nabla :  {\mathcal A} ^{p,q} \rightarrow  {\mathcal A} ^{p+1,q},$$
which give a skew-derivation on the algebra ${\mathcal A} ^{*,*}$.  To simplify notation, we often suppress writing the wedge.

When we apply $\nabla $ to a smooth section $V$ of $E$, it gives the {\em covariant differential\/} $\nabla V \in {\mathcal A} ^{1,1}$, which is defined in agreement with (\ref{E:covdiff}) by
$$\nabla V = \sum_{i=1}^n \left(dv_i + \sum_{j=1}^n \omega _{ij} v_j \right) e_i,$$
where $V$ is expressed in terms of a moving orthonormal frame $e_1, \ldots , e_n$ defined on an open subset $U\subset M$ by
$$V = \sum_{i=1}^nv_i e_i,$$
the $\omega_{ij}$'s being the connection one-forms determined by the moving frame.  The curvature forms of the connection,
$$\Omega _{ij} = d\omega _{ij} + \omega _{ik} \wedge \omega _{kj}$$
fit together into a {\em curvature operator\/}
$${\cal R}  = - \frac{1}{4} \sum _{i,j = 1}^n \Omega _{ij} e_i \wedge e_j \in \Omega^{2,2},$$
which is independent of the choice of moving frame, like $\nabla V$.  All the summands in the expression
$$\Phi (V,\nabla ) = - |V|^2 + \nabla V + {\cal R} \in {\mathcal A}^{0,0} \oplus {\mathcal A} ^{1,1} \oplus {\mathcal A} ^{2,2}$$
have even total degree and therefore commute with each other.  Thus all terms in the expanded power series
$$\hbox{exp}(\Phi (V,\nabla )) = I + [- |V|^2 + \nabla V + {\cal R}] + (1/2)[- |V|^2 + \nabla V + {\cal R}]^2 + \cdots $$
commute, and we can write
$$\hbox{exp}(\Phi (V,\nabla )) = e^{-|V|^2}[I + (\nabla V)+ (1/2) (\nabla V)^2 + \cdots ] [ I + {\cal R} + (1/2) ({\cal R})^2 + \cdots ].$$
The two infinite series within brackets have only finitely many terms because
${\mathcal A} ^{p,q} = 0$ when $p > \dim (M)$ or $q > n = \hbox{rank}(E)$.

The orientation of $E$ determines a volume element $\star 1 \in \Lambda ^n (E)$, which can be expressed in terms of a positively oriented moving orthonormal frame as
$$\star 1 = e_1 \wedge \cdots \wedge e_n.$$
Following Quillen \cite{Q}, \cite{MQ}, we define the {\em supertrace\/}
$$\hbox{Tr}^M_s: {\mathcal A}^p(M,\Lambda ^*E) \rightarrow {\mathcal A}^p(M),$$
the space of ordinary $p$-forms on $M$, by projection on the $\star 1$-component; thus
$$\hbox{Tr}^M_s\left( \sum _{i_1 < \cdots < i_k} \phi _{i_1\cdots i_k} e_{i_1} \wedge \cdots \wedge e_{i_k}\right) = \phi _{1\cdots n}.$$
For any choice of section $V$, we can construct the $n$-form
\begin{equation}\tau (V , \nabla ) = \frac{(-1)^{[n/2]}}{\pi^{n/2}} \hbox{Tr}^M_s [\hbox{exp}(\Phi (V,\nabla ))] \label{eq:tauv}\end{equation}
on $M$, where $[n/2]$ represents the largest integer $\leq n/2$.  This form can be constructed not only for the bundle $E$ itself, but also for the pullback bundle via the bundle projection $\pi : E \rightarrow M$.  This pullback bundle
$$\pi^*E = \{ (e_1,e_2) \in E \times E : \pi (e_1) = \pi (e_2) \}$$
inherits an inner product, an orientation, and a pullback connection $\pi ^*\nabla $ from $E$, and it possesses a \lq\lq tautological" section
$$T: E \rightarrow \pi ^* E, \qquad T(e) = (e,e).$$

\vskip .1 in
\noindent
{\bf Definition.}  The {\em Thom form\/} on $E$ is the differential form
\begin{equation} \tau (\nabla ) = \frac{(-1)^{[n/2]}}{\pi^{n/2}} \hbox{Tr}^M_s [\hbox{exp}(\Phi (T,\pi ^*\nabla ))].\label{E:thomform}\end{equation}

\vskip .1 in
\noindent
This Thom form $\tau (\nabla )$ depends upon the choice of metric connection $\nabla $, which is expressed in terms of the connection forms $\omega _{ij}$ with respect to a moving frame, but $\pi ^*\omega _{ij}$ and $\pi ^*\Omega _{ij}$ restrict to zero on each fiber.  Hence if we represent a general point in the fiber as $\sum t_i e_i$, the restriction of $\tau (\nabla )$ to the fiber is
$$\frac{(-1)^{[n/2]}}{\pi^{n/2}} e^{-(t_1^2 + \cdots + t_n^2)} \hbox{Tr}^M_s\left( \frac{1}{n!} \left( \sum dt_i e_i\right)^{n} \right),$$
and by a short calculation, one verifies that this is just
$$\frac{1}{\pi^{n/2}} e^{-(t_1^2 + \cdots + t_n^2)} dt_1 \wedge \cdots \wedge dt_{n},$$
since we must make $n(n-1)/2$ changes of sign when permuting $dt_i$'s and $e_j$'s and $(-1)^{n(n-1)/2} = (-1)^{[n/2]}$.  Thus $\tau (\nabla )$ integrates to one on each fiber, is rapidly decreasing in the fiber direction, and specializes to the Thom form constructed before for the case of the trivial bundle with trivial connection.  Moreover, if $V$ is a section of $E$, then $V^*(\tau (\nabla ))$ is the form $\tau (V , \nabla )$ defined by (\ref{eq:tauv}).

\vskip .1 in
\noindent
{\bf Lemma 4.1.}\begin{sl}  The Thom form on $E$ is closed.  \end{sl}

\vskip .1 in
\noindent
Proof:  If $V$ is a section of $\pi ^*E$, let $\iota (V) : {\mathcal A}^{p,q} \rightarrow {\mathcal A}^{p,q-1}$ be the interior product, the skew-derivation such that $\iota (V)(W) = \langle V, W \rangle$ when $W$ is a section of $\pi ^*E$.  It is then relatively straightforward to verify that
\begin{equation} \nabla (-|T|^2) = 2 \iota (T)\nabla T, \qquad \nabla (\nabla T) = 2 \iota (T){\cal R}, \qquad \nabla {\cal R} = 0, \label{E:lemma} \end{equation}
the last equality following from the Bianchi identity
$$d\Omega _{ij} = \sum \Omega _{ik} \wedge \omega _{kj} - \sum \omega _{ik} \wedge \Omega _{kj}.$$
It follows from (\ref{E:lemma}) that
$$\nabla \Phi (T,\nabla ) = 2 \iota (T) \Phi (T,\nabla ),$$
and since $\nabla $ and $\iota (T)$ are both skew-derivations, that
$$\nabla (\hbox{exp}(\Phi (T,\nabla ))) = 2 \iota (T) \hbox{exp}(\Phi (T,\nabla )).$$
Now one applies the identities
$$\hbox{Tr}^M_s \circ \iota (T) = 0, \qquad \hbox{Tr}^M_s \circ \nabla = d \circ \hbox{Tr}^M_s$$
to conclude that the Thom form $\tau (\nabla )$ on $E$ is indeed closed.  QED

\vskip .1 in
\noindent
One can check that the cochain complex of smooth forms on $E$ which are rapidly decreasing in the fiber direction is cochain homotopy equivalent to the cochain complex of smooth forms on $E$ which have compact support in the fiber direction.  Hence we can regard the Thom form as repreenting a cohomology class in $H^*_{cv}(E;{\mathbb R})$, the de Rham cohomology of $E$ with compact support in the fiber or \lq\lq vertical" direction, as described in Bott and Tu \cite{BT}, page~61, and indeed using the Mayer-Vietoris technique described in \cite{BT}, one shows that wedging with the Thom form yields an isomorphism
$$\cdot \cup [\tau (\nabla )] : H^k(M;{\mathbb R}) \longrightarrow H^{k+n}_{cv}(E;{\mathbb R}),$$
which is called the {\em Thom isomorphism\/}.  

If $V : M \rightarrow E$ is a smooth section, we recover the form $\tau (V,\nabla )$ defined earlier by (\ref{eq:tauv}) by pulling back via $V$: $\tau (V,\nabla ) = V^*(\tau (\nabla ))$.  In particular, if $n$ is even, say $n = 2m$,
$$\tau(0,\nabla ) = \frac{(-1)^m}{m!\pi^m}\hbox{Tr}^M_s({\cal R}^m) = \frac{1}{(2\pi )^m} \hbox{Pf}(\Omega ),$$
the last equality following from (\ref{E:pfaff}) and the definition of ${\cal R}$, while if $n$ is odd, $\tau(0,\nabla ) = 0$ since $\hbox{Tr}^M_s({\cal R}^m)$ vanishes.  In either case, we get the Gauss-Bonnet integrand.

To prove the Gauss-Bonnet Theorem for closed manifolds, we suppose that $V$ is a section of $E$ such that at each zero $p$ of $X$,
$$\nabla V(p): T_pM \rightarrow E_p$$
is an isomorphism (in which case we can say that the zeros of $V$ are {\em nondegenerate\/}).  The {\em rotation index\/} of $\nabla V$ at such a zero $p$ is defined to be either $1$ or $-1$ depending upon whether $\nabla V(p)$ preserves or reverses orientation.  A short calculation shows that as $s \rightarrow \infty$, $(sV)^* (\tau _E)$ becomes concentrated near the zeros of $V$, and in fact,
\begin{equation}\int _M \tau (sV, \nabla ) \rightarrow \sum (\hbox{rotation indices of $V$ at its zeros}).\label{E:rotationindex}\end{equation}

This calculation can be carried out in terms of normal coordinates $(x_1, \ldots x_n)$ centered at a zero $p$ for $V$ and defined in a small neighborhood $B_\epsilon (p)$ of radius $\epsilon $ about $p$.  If we write
$$V = \sum a_{ij}x_i\frac{\partial }{\partial x_j} + (\hbox{higher order terms}),$$
then since the $\omega _{ij}$'s vanish at $p$,
$$\nabla V = \sum a_{ij} dx_i e_j + (\hbox{higher order terms}), \quad \hbox{where} \quad e_j = \left. \frac{\partial }{\partial x_j}\right|_p.$$
Thus
\begin{multline*} (\nabla V)^n = \sum a_{i_1j_1} dx_{i_1} e_{j_1} \cdots a_{i_nj_n} dx_{i_n} e_{j_n} \\ = (-1)^{[n/2]}\sum a_{i_1j_1} \cdots a_{i_nj_n} (dx_{i_1} \wedge \cdots \wedge dx_{i_1})(e_{j_1} \wedge \cdots \wedge e_{j_n}),\end{multline*}
the sign coming from commuting the $dx_i$'s and the $e_j$'s.  Using the definition of determinant, we conclude that
\begin{multline*} \frac{1}{n!}\left( \nabla V \right)^{n} = (-1)^{[n/2]} \det (a_{ij})(dx_1 \wedge \cdots dx_n) (e_1 \wedge \cdots \wedge e_n) \\ + (\hbox{higher order terms}),\end{multline*}
Hence
\begin{multline*} \hbox{Tr}^M_s\left( e^{-s^2 |V|^2} \frac{s^{n}}{n!}\left( \nabla V \right)^n\right) = (-1)^{[n/2]} s^{n} e^{-s^2 |V|^2} \det (a_{ij})(dx_1 \wedge \cdots dx_n) \\
+ (\hbox{higher order terms}).\end{multline*}
After the change of variables $y_i = \sum sa_{ij} x_j$, this becomes
$$(-1)^{[n/2]} e^{-(y_1^2 + \cdots + y_n^2)} (dy_1 \wedge \cdots dy_n) + (\hbox{higher order terms}).$$
Note that the ball $B_\epsilon (p)$ of radius $\epsilon $ gets rescaled to a ball of radius $s \epsilon $ in the $(y_1, \ldots y_n)$-coordinates.  We find that when $s$ is large, 
\begin{multline*} \int _{B_\epsilon (p)}\tau (sV, \nabla ) = \frac{1}{\pi^{n/2}} \int _{B_\epsilon (p)} e^{-(y_1^2 + \cdots + y_n^2)} (dy_1 \wedge \cdots dy_n)\\
+ (\hbox{higher order terms}).\end{multline*}
As $s \rightarrow \infty$ (and hence $s \epsilon \rightarrow \infty$), we can use (\ref{E:gaussian}) to obtain
$$\hbox{lim}_{s \rightarrow \infty} \int_{B_\epsilon (p)} \tau (sV, \nabla ) = \begin{cases} 1, & \hbox{if $\det (a_{ij}) > 0$,} \cr -1, & \hbox{if $\det (a_{ij}) < 0$,} \end{cases} $$
and this is just the rotation index of $V$ at $p$.  Adding together the contributions at all the zeros of $V$ yields (\ref{E:rotationindex}).

On the other hand, as $s \rightarrow 0$, $\tau (sV, \nabla)$ approaches the Gauss-Bonnet integrand when $n$ is even or zero if $n$ is odd.  The integral is constant as a function of $s$, because the integrands are pullbacks of a closed form under a smooth homotopy of maps, and hence
\begin{equation} \frac{(-1)^{[n/2]}}{\pi^{n/2}} \int _M \hbox{Tr}^M_s(e^{\cal R}) = \sum (\hbox{rotation indices of $V$ at its zeros}). \label{E:poincarehopf} \end{equation}
The left-hand side is independent of the choice of section $V$ while the right-hand side is independent of the connection $\nabla $.  Thus neither side depends on either the section $V$ or the connection $\nabla $ and both sides must give an expression for a topological invariant.  If $E = TM$, we can take $V$ to be the gradient of a Morse function on $M$, which identifies this invariant with the Euler characteristic $\chi (M)$ by standard Morse theory, as presented for example in Milnor \cite{Mil}, finishing the proof of the generalized Gauss-Bonnet Theorem for closed manifolds.

\section{Gauss-Bonnet for manifolds with boundary}
\label{S:integrationoverfibers}

The Hazzidakis formulae needed for the proof of the Main Theorem are obtained by generalizing the Gauss-Bonnet Theorem first to compact smooth manifolds with boundary, and then to Riemannian polyhedra having corners, such as the polyhedra $P_a$ defined by (\ref{E:Riemannianpolyhedron}).

To extend our previous argument to manifolds with boundary, we consider the unit sphere bundle whose total space is
$$E^1 = \{ v \in E: \langle v, v \rangle = 1 \}.$$
Corresponding to the Euler form $e(\nabla )$ defined for a given metric connection $\nabla $ on the oriented bundle $E$ by (\ref{E:eulerform}) on $M$, there is a {\em transgressed Euler form\/} on $E^1$, an $(n-1)$-form $Te(\nabla )$ on $E^1$ such that
\begin{equation} d(Te(\nabla )) = \pi^* e(\nabla), \quad \hbox{where} \quad \pi : E^1 \longrightarrow M \label{E:transgression} \end{equation}
is the projection.\footnote{It is analogous to the transgressed Pontrjagin form constructed by Chern and Simons \cite{ChS} in their treatment of Chern-Simons invariants.}  This transgressed Euler form is constructed by integrating the Thom form of \S \ref{S:gaussbonnet} over fibers of $E$, each fiber being a ray $\{ tv : t \in [0,\infty) \}$ generated by an element $v \in E^1$.

\vskip .1 in
\noindent
{\bf Integration over fibers.}  Integration over fibers is discussed in pages 61-65 of \cite{BT}, with an extension to the case in which each fiber has a nontrivial boundary in Theorem 3.10 in \cite{B}:  Let $E$ be a smooth fiber bundle over the base $M$ with fiber a smooth manifold $F$ of dimension $q$, which may possess a smooth boundary $\partial F$.  Suppose that  $\Theta $ is a smooth $(p+q)$-form on $E$ with compact support in the fiber direction, and that the restriction of $\Theta $ to $\partial F$ also has compact support in the fiber direction.  (More generally, we can allow forms which are rapidly decreasing in the fiber direction.) We can then integrate $\Theta $ over $F$ to obtain a $p$-form $\Theta /F$ on $M$.  In terms of local coordinate $(x_1, \ldots , x_s,t_1, \ldots , t_q)$, where $(x_1, \ldots , x_s)$ are coordinates on the base and $(t_1, \ldots t_q)$ restrict to coordinates on the fiber,
\begin{multline*} \Theta = f(t,x) dx_{i_1} \wedge \cdots \wedge dx_{i_p} \wedge dt_1 \wedge \cdots \wedge dt_q \\ \mapsto \quad \Theta/F = \left( \int_F f(t,p) dt_1 \ldots dt_q \right)dx_{i_1} \wedge \cdots \wedge dx_{i_p},\end{multline*}
while forms that do not contain the factor $dt_1 \wedge \cdots \wedge dt_q$ integrate to zero.  Moreover, we easily verify that
\begin{multline} (d\Theta )/F = d(\Theta /F) + (-1)^p \Theta /\partial F, \qquad \hbox{so that} \\ d\Theta = 0 \quad \Rightarrow \quad d(\Theta /F) = - (-1)^p \Theta /\partial F. \label{E:fibint} \end{multline}
We use the slant product notation and call $\Theta /F$ the {\em slant product\/} of $\Theta $ with $F$ in agreement with the slant product on the chain and cochain level in algebraic topology, as described in Spanier \cite{Sp}, Chapter 6, \S 1.

In the case where $E$ is a vector bundle of rank $n$ over $M$ with fiber ${\mathbb R}^n$, the Mayer-Vietoris technique of \cite{BT} shows that integration over the fiber induces an isomorphism $H^{p+n}_{cv}(E;{\mathbb R}) \rightarrow H^p(M;{\mathbb R})$ which is an inverse to the Thom isomorphism described before.

\vskip .1 in
\noindent
{\bf The transgressed Euler form.}  To construct the transgressed Euler form, a smooth $(n-1)$-form on $E^1$, we pull the bundle $E$ back to a bundle $\pi ^*E$ over $E^1$ via the projection
$$\pi : E^1 \longrightarrow M.$$
As in the argument from \S \ref{S:gaussbonnet}, we then have a tautological section
$$T : E^1 \longrightarrow \pi^*E, \quad \hbox{defined by} \quad T(e) = (e,e).$$
For each $e \in E^1$, we consider the ray
$$F_e = \{ t T(e) = (e,te) \in (\pi^*E)_e : t \in [0,\infty ) \}.$$
As $e$ ranges over $E^1$, the family of rays $e \mapsto F_e$ generate a fiber bundle over $E^1$ with fibers diffeomorphic to ${\mathbb R}^{\geq 0} = \{ t \in {\mathbb R} : t \geq 0 \}$.  When applying integration over the fiber we take $p = 2n-1$, so (\ref{E:fibint}) yields
$$d(\Theta /F) = \Theta /\partial F = - \Theta|\partial F = - (\hbox{the restriction of $\Theta$ to $\partial F$}),$$
the minus occurring because with the standard orientation of $[0,\infty)$, zero is the lower limit of integration.  The {\em transgressed Euler form\/} on $E^1|\partial M$ is
$$Te(\nabla ) = - \pi^*\tau(\nabla )/F,$$
where $\tau(\nabla )$ is the Thom form.  Note that it satisfies the key identity (\ref{E:transgression}).

\vskip .1 in
\noindent
{\bf The geodesic curvature form.} If $V$ is a unit-length section of $E^1|\partial M$, we call
\begin{equation} \Pi (\nabla ,V) = - V^*Te(\nabla ). \label{E:geodcurvform} \end{equation}
the {\em geodesic curvature form\/} corresponding to $V$ along $\partial M$.  It follows from (\ref{E:transgression}) that the restriction of $dTe(\nabla )$ to $E^1|\partial M$ is zero, and hence
$$\int _{\partial M} \Pi (\nabla ,V) = - \int _{\partial M}V^*Te(\nabla )$$
depends only on the homotopy class of the unit-length section $V : \partial M \longrightarrow E^1|\partial M$.  We claim that this geodesic curvature form fits into  a generalization of the Lccal Gauss-Bonnet Theorem described on page~272 of do Carmo \cite{dC}:

\vskip .1 in
\noindent
{\bf Theorem 5.1. (Gauss-Bonnet for Manifolds with Boundary).} \begin{sl} Suppose that $M$ is a smooth compact oriented Riemannian manifold with smooth boundary $\partial M$ and that $V$ is a section of $TM$ which has finitely many nondegenerate zeros within $M - \partial M$, has unit length along $\partial M$ and is homotopic to the outward-pointing unit normal $N$ along $\partial M$.  Moreover, suppose that $\nabla $ is a metric connection on $TM$ with curvature $\Omega $ and $\chi (M)$ is the Euler characteristic of $M$.
\begin{enumerate}
\item If $M$ is even-dimensional, say $\dim M = 2m$, then
\begin{equation} \frac{1}{(2\pi)^m} \int_{M} \hbox{Pf}(\Omega ) + \int_{\partial M} \Pi (\nabla ,V) = \chi (M). \label{E:gbwithboundarye}\end{equation}
\item If $M$ is odd-dimensional, then
$$ \int_{\partial M} \Pi (\nabla ,V) = \chi (M). $$
\end{enumerate}
\end{sl}

\vskip .1 in
\noindent
This was proven by Chern (see (19) on page 678 of \cite{Ch45}), but note that Chern takes $V$ to equal the inward-pointing unit normal $N$, while we require that $V$ be only homotopic to the outward-pointing unit normal $N$.  Since we take $N$ to be outward-pointing, our signs differ from those in some of Chern's formulae.

\vskip .1 in
\noindent
{\bf Proof of Theorem~5.1:}  Let $V$ be a smooth vector field on $M$ with finitely many nondegenerate zeros in $M - \partial M$ such that $V$ has unit length along $\partial M$.  If $s$ is a large element of $[0,\infty)$, we let
$$F_s = \{ t V(p) : p \in M, t \in [0,s] \}.$$
Then the boundary of $F_s$ consists of three pieces, the zero section $Z \subseteq TM$, the image $\sigma_s(M)$ of the section
$$\sigma _s : M \longrightarrow E, \quad \hbox{defined by} \quad \sigma _s(p) = s V(p),$$
and the part $F_s|\partial M$ of $F_s$ over $\partial M$.  Indeed, with appropriate orientations on these three pieces, we have
$$\partial F_s = \sigma_s(M) + F_s|\partial M - Z.$$
The Thom form $\tau (\nabla )$ is closed, so it integrates to zero over the boundary $\partial F_s$ of the contractible polyhedron $F_s$, and after shifting the integral over $Z$ to the other side, we obtain
$$\int_Z \tau (\nabla) = \int _{F_s|\partial M} \tau (\nabla ) + \int _M \sigma _s^*\tau(\nabla ).$$
The Gauss-Bonnet argument of \S \ref{S:gaussbonnet} shows that as $s \rightarrow \infty$ the three terms in this expression converge to the corresponding terms of
\begin{equation}  \frac{1}{(2\pi)^m} \int_{M} \hbox{Pf}(\Omega ) = - \int_{\partial M} \Pi (\nabla ,V) + \sum (\hbox{rotation indices of $V$ at zeros}). \label{E:poincarehopf1} \end{equation} 

The first term in (\ref{E:poincarehopf1}) depends only on the Riemannian geometry on $M$, not on the choice of $V$.  The second term on the right depends only on the singularities of $V$ within $M - \partial M$, not on the Riemannian metric.  Since $V$ is homotopic to the outward-pointing unit normal within a tubular neighborhood of $\partial M$ which avoids the singularities, we can replace $V$ with a new vector field with the same singularities as $V$ that also satisfies the condition that it restrict to the outward-pointing unit normal along $\partial M$.

Using a tubular neighborhoods of $\partial M$ within $M$, we can arrange that $V$ is the gradient of $f$, where $f$ is a smooth Morse function on $M$ with $f \leq 1$ on $M$ and $f = 1$ along $\partial M$.  With this choice of $V$, Morse theory implies that the second term on the right side of (\ref{E:poincarehopf1}) must equal the Euler characteristic of $M$, as in the proof of the Gauss-Bonnet formula for closed manifolds.  QED

\vskip .1 in
\noindent
{\bf Remark~5.2. An explicit formula for the geodesic curvature form over $\partial M$ when $n$ is even:}  To apply Theorem~5.1 we often need a formula for the geodesic curvature form $\Pi(\nabla , V)$ along $\partial M$.  We can use (\ref{E:thomform}) to calculate such an expression, and focus first on the case where $n$ is even.  Since
$$\nabla (tT) = (dt) T + t \nabla T,$$
and $dt T$ commutes with every term of $\hbox{exp}(\Phi (T,\pi ^*\omega ))$, when $T$ denotes the tautological section of $\pi ^*E$, an integration of the Thom form $\pi^*\tau (\nabla )$ over the fiber $F$ yields the following expression for the geodesic curvature form of (\ref{E:geodcurvform}):
\begin{equation} \Pi(\nabla , V) = \frac{(-1)^{[(n-1)/2]}}{\pi^{(n/2)}} \int _0^\infty e^{-t^2} dt \hbox{Tr}^M_s[V \exp (t\nabla V) \exp ({\mathcal R}^M )].\label{E:fiberint} \end{equation}
Note that when $n$ is even the integrand in (\ref{E:fiberint}) is odd in $t$, and we can rewrite this integral as
$$\frac{1}{2} \frac{(-1)^{[(n-1)/2]}}{\pi^{(n/2)}} \left[ \int _0^\infty - \int _{-\infty}^0 \right] e^{-t^2} dt \hbox{Tr}^M_s[V \exp (t\nabla V) \exp ({\mathcal R}^M )].$$
This implies an important property of geodesic curvatures forms:
\begin{equation}\hbox{$n$ even} \quad \Rightarrow \quad \Pi(\nabla ,-V) = \Pi(\nabla , V).\label{E:symneven}\end{equation}
In other words, when $n$ is even, $\Pi(\nabla , V)$ is invariant under the reflection in the fiber of the line bundle
$$F = \{ t V(p) : p \in \partial M, t \in [{\mathbb R} \} \longrightarrow \partial M$$
which sends $V$ to $-V$.

\vskip .1 in
\noindent
{\bf Remark~5.3. An explicit formula for the geodesic curvature form over $\partial M$ when $n$ is odd:}  When $n$ is odd, the integrand in (\ref{E:fiberint}) is even in $t$, so we can write
$$\Pi(\nabla , V) = \frac{1}{2} \frac{(-1)^{[(n-1)/2]}}{\pi^{(n/2)}} \int _{-\infty}^\infty e^{-t^2} dt \hbox{Tr}^M_s[V \exp (t\nabla V) \exp ({\mathcal R}^M )],$$
and only terms in $\exp (t\nabla V)$ which are even in $t$ can survive.  We can drop the $N$ within the supertrace if we replace the supertrace $\hbox{Tr}^M_s$ for $M$ by the supertrace $\hbox{Tr}^{\partial M}_s$ for $\partial M$ and multiply by $(-1)^n$.  The integral then becomes
$$(-1)^n \sum _{k=0}^\infty \left[\int _{-\infty}^\infty  t^{2k} e^{-t^2} dt \right] \hbox{Tr}^{\partial M}_s \left[ \frac{(\nabla V))^{2k}}{(2k)!} \exp ({\mathcal R}^M) \right],$$
where we understand that the curvature operator ${\mathcal R}^M$ has been restricted to the exterior algebra on the submanifold $\partial M$.  But looking up the integral in an appropriate table\footnote{See for example \cite{AS}, page 302, formula 7.4.4.} yields
$$ \int _{-\infty}^\infty  t^{2k} e^{-t^2} dt = \frac{\sqrt \pi}{2^{k}}(2k-1)(2k-3) \cdots 3 \cdot 1,$$
so we obtain the formula
\begin{equation} \Pi(\nabla ,V)  = (-1)^n \sqrt \pi \hbox{Tr}^{\partial M}_s \left[ \exp \left( \frac{(\nabla V)^2}{4} +{\mathcal R}^M \right) \right]. \label{E:geodesiccurvforboundaryMmod}\end{equation}

In the case where $V$ is the outward pointing unit normal $N$, this formula becomes
\begin{equation} \Pi(\nabla ,N)  = (-1)^n \sqrt \pi \hbox{Tr}^{\partial M}_s \left[ \exp \left( \frac{(\nabla N)^2}{4} +{\mathcal R}^M \right) \right]. \label{E:geodesiccurvforboundaryM}\end{equation}
We can analyze this in terms of an adapted moving frame $(e_1, \ldots e_{n-1}, e_n)$ along $\partial M$ with $(e_1, \ldots e_{n-1})$ tangent to $\partial M$ and $e_n = N$.  The curvature forms $\Omega ^M_{ij}$ for any connectionn on $M$ are related to the curvature forms $\Omega ^{\partial M}_{ab}$ for the induced metric connection on $\partial M$, where $1 \leq a,b \leq n-1$, by the Gauss equation
\begin{equation} \Omega ^{\partial M}_{ab} = \Omega ^M_{ab} + \omega _{na} \wedge \omega _{nb}, \quad \hbox{for $1 \leq a,b \leq n-1$.} \label{E:Gaussequation} \end{equation}
Note that this Gauss equation holds for arbitrary metric connections not just the Levi-Civita connection which has vanishing torsion.  The curvature operators ${\mathcal R}^M$ and ${\mathcal R}^{\partial M}$ are defined by
$${\mathcal R}^M = - \frac{1}{4}\sum _{i,j} \Omega ^M_{ij} e_i \wedge e_j, \quad {\mathcal R}^{\partial M} = - \frac{1}{4}\sum _{a,b} \Omega ^{\partial M}_{ab} e_a \wedge e_b,$$
while
$$\frac{(\nabla N)^2}{4} = \frac{1}{4} \sum _{a,b} \omega _{na} e_a \wedge \omega _{nb} e_b = - \frac{1}{4} \sum _{a,b} (\omega _{na} \wedge \omega _{nb}) (e_a  \wedge e_b),$$
so the Gauss equation (\ref{E:Gaussequation}) translates to
\begin{equation} {\mathcal R}^{\partial M} = {\mathcal R}^M + \frac{1}{4} (\nabla N)^2. \label{E:nsquare}\end{equation}
Substituting this into (\ref{E:geodesiccurvforboundaryM}) yields one-half the geodesic curvature form for $\partial M$.  What we have proven is:

\vskip .1 in
\noindent
{\bf Lemma 5.4.}\begin{sl}  When $M$ has odd dimension and we choose $V$ to be the outward-pointing unit normal $N$, the geodesic curvature form for $\partial M$ is simply one-half the Gauss-Bonnet integrand expressed in terms of the metric connection induced on $\partial M$. \end{sl}

\vskip .1 in
\noindent
When applied to the Levi-Civita connection, this Lemma is a special case of the \lq\lq new integral formula" treated in \S 3 of \cite{Ch45}.  Indeed, Chern observes that this Lemma together with the Gauss-Bonnet formula applied to $\partial M$ gives the following Corollary of the second assertion in Theorem~5.1: If $M$ is a smooth odd-dimensional compact oriented Riemannian manifold with smooth boundary $\partial M$, then
\begin{equation} \chi (M) = \frac{1}{2}{ \chi (\partial M}). \label{E:gbbdryoddcase}\end{equation}

In our later application of the geodesic curvature form, it will be convenient to use (\ref{E:geodesiccurvforboundaryMmod}) in the case where the unit-length vector field $V$ is only homotopic to the outward pointing unit normal $N$.  Note that in the general case, the formula (\ref{E:geodesiccurvforboundaryMmod}) shows that in contrast to (\ref{E:symneven}), we have
\begin{equation}\hbox{$n$ odd} \quad \Rightarrow \quad \Pi(\nabla ,-V) = - \Pi(\nabla , V).\label{E:symnodd}\end{equation}
In other words, when $n$ is odd, $\Pi(\nabla , V)$ is anti-invariant under reflection in the fiber of the line bundle
$$F = \{ t V(p) : p \in \partial M, t \in [{\mathbb R} \} \longrightarrow \partial M$$
which sends $V$ to $-V$.

\vskip .1 in
\noindent
{\bf Remark~5.5.  Hazzidakis for surfaces.}  We can derive the Hazzidakis formula in the special case of surfaces as an application of Theorem~5.1 with suitable adjustments at the corners.  For the principal coordinates $(y_1, y_2)$, we have the coordinate vector fields
$$Y_1= \frac{\partial }{\partial y_1} \quad \hbox{and} \quad Y_2= \frac{\partial }{\partial y_2},$$
while the corresponding unit-length asymptotic vector fields are 
$$Z = Y_1 + Y_2 \quad \hbox{and} \quad W = Y_1 - Y_2.$$
It then follows from (\ref{eq:hilb28}) and (\ref{eq:hilb29}) that
$$\omega _{12}(W) = \phi _{12}(W),$$
which translates into the condition that
$$\nabla _WZ = \widehat \nabla _WZ,$$
where $\nabla $ is the Levi-CIvita connection and $\widehat \nabla $ is the connection in $NM$ pulled back via $h \circ f$.  But $Z$ is parallel with respect to $\widehat \nabla $ as explained at the end of \S \ref{S:adaptedframe}, so
$$\nabla _WZ = 0 \quad \hbox{and similarly} \quad \nabla _ZW = 0.$$

We apply these formulae to the coordinate rectangle $R$ of (\ref{E:asymptoticrectangle}).  Along an edge in the boundary of $R$ which is tangent to $Z$, we let $V = \pm W$, choosing the sign so that $V$ is outward-pointing.  Similarly, along an edge which is tangent to $W$, we let $V = \pm Z$, choosing the sign so that $V$ is outward-pointing.  Then the integrals along edges in the Gauss-Bonnet formula (\ref{E:gbwithboundarye}) for $R$ vanish, but we need to add correction terms at each vertex, where the vector field along the boundary jumps.\footnote{This correction is easily calculated for a surface; in the higher-dimensional cases we extend the process of integration over the fiber to deal with multi-dimensional corners.}  The correction needed at a vertex when we move along the boundary in a counterclockwise direction  is
$$- \frac{\varepsilon _i}{2\pi}, \quad \hbox{where $\varepsilon _i$ is the exterior angle at the $i$-th vertex,}$$
with $i$ running from $1$ to $4$.  We set $K = -1$ and $\chi (M) = 1$, so if $(z,w)$ are the coordinates corresponding to $Z$ and $W$,
$$\hbox{Pf}(\Omega ) = - \sin \theta dz dw,$$
and (\ref{E:gbwithboundarye}) becomes
$$- \frac{1}{2\pi} \int_R \sin \theta dzdw = - \sum_{i=1}^4 \frac{\varepsilon _i}{2\pi} + 1,$$
which is just (\ref{E:hilb12}) after a change of sign.

This argument suggests a strategy for proving the Main Theorem, namely showing that covariant derivative terms in the Allendoerfer-Weil version of the Gauss-Bonnet formula for a Riemannian polyhedron vanish.  We follow a modified version of  that strategy in the following sections.

\vskip .1 in
\noindent
{\bf Remark~5.6. The special case of an $n$-disk:}. Of special importance to us is the case in which $M$ is the $n$-disk $D^n$ with boundary $S^{n-1}$ because this is the case modified by Allendoerfer and Weil to treat $n$-cells with corners on the boundary.  In this case the geodesic curvature form $\Pi (\nabla ,V)$ along $\partial D^n$ represents the positively oriented generator of $H_{n-1}(S^{n-1};{\mathbb Z})$.  Indeed, since $\chi (D^n) = 1$, the Gauss-Bonnet formula of Theorem~5.1 reduces to
\begin{equation} \frac{1}{(2\pi)^m} \int_{D^n} \hbox{Pf}(\Omega ) + \int _{S^{n-1}} \Pi (\nabla ,V) = 1 \label{E:GBDeven} \end{equation}
when $n$ is even, or
\begin{equation} \int _{S^{n-1}} \Pi (\nabla ,V) = 1 \label{E:GBDodd} \end{equation}
when $n$ is odd.  Here $V$ is any unit-length vector field along $\partial D^n =S^{n-1}$ which is homotopic to the outward-pointing unit normal.  Such a vector field along $\partial D^n$ could be extended to a smooth vector field on $D^n$ with a single nondegenerate zero at the center of $D^n$, the zero having rotation index one.

\section{Dual CW decompositions for $S^{n-1}$}
\label{S:dualCW}

Our next goal is to extend the Gauss-Bonnet formulae of (\ref{E:GBDeven}) and (\ref{E:GBDodd}) so that they apply to the asymptotic polyhedra of (\ref{E:Riemannianpolyhedron}),
$$P_a = \{ p \in {\mathbb H}^n : |\pm y_1(p)\pm \cdots \pm y_n(p) | \leq a, \hbox{ for all choices of signs}\},$$
in which $(y_1, \ldots , y_n)$ is a preferred set of principal coordinates on a contractible open subset $U$ of our hyperbolic manifold $M$, which is isometrically imbedded in ${\mathbb E}^{2n-1}$.  The boundary of $P_a$ has an obvious triangulation in which the zero-dimensional simplices are the $2n$ vertices described in terms of the principal coordinates as
$${\bf a}_1 = (a,0, \ldots ,0), \qquad \ldots, \qquad {\bf a}_n = (0,0, \ldots ,a),$$
$$-{\bf a}_1 = (-a,0, \ldots ,0), \qquad \ldots, \qquad -{\bf a}_n = (0,0, \ldots ,-a).$$
We will eventually choose $V$ along each $(n-1)$-dimensional face of this triangulation to be a unit-length asymptototic vector field transverse to this face and pointing outward.  But this requires modifying $V$ along lower-dimensional simplices where these choices disagree.  Thus we need to adjust the Gauss-Bonnet formula of Remark~5.6 to account for changes of $\Pi (\nabla ,V)$ along simplices of dimension $\leq n-2$.  We will see that this can be done by further integration over the fiber, extending the process used in the previous section.

The  principal frame $(e_1, \ldots e_n)$ over $U$ establishes a trivialization $TM|U \cong U \times {\mathbb R}^n$, and a product structure on the corresponding unit tangent bundle $T^1M|U \cong U \times S^{n-1}$.   This trivialization enables us to regard the triangulation of the boundary $\partial P_a$ of an asymptotic polyhedron $P_a$ as having a dual CW decomposition within the fiber $T^1_{p_0}M$ where $p_0$ is the center of the polyhedron, with coordinates
$$y_1 = y_2 = \cdots = y_n = 0.$$
We remark that a similar construction could be carried out for any Riemannian polyhedron homeomorphic to an $n$-disk and lying within a convex open subset of a preferred coordinate system.  The advantage to the asymptotic polyhedra we have chosen is that they have a large group $G$ of symmetries which induce an isometric action on the fibers of the unit tangent bundle $T^1M$.

Recall the definition of our symmetry group $G$; it is the abelian group of order $2^n$ generated by the coordinate reflections $g_i$, for $1 \leq i \leq n$, which are defined by
$$g_i(y_j) = (1 - 2 \delta _{ij}) y_j.$$
Each $g  \in G$ corresponds to a map
$$g : \{ 1 , \ldots ,n \} \rightarrow \{ -1,1 \} \qquad  \hbox{such that} \qquad g(y_i) = g(i) y_i.$$
It is important to note that the elements of $G$ do not necessarily act as isometries of the asymptotic polyhedron $P_a$ itself, but $G$ does act as isometries on $TM|P_a$, and  we will see that this induces an action on the space of metric connections on $TM|P_a$, an action that will be crucial to our proof.

We let ${\cal I}_k$ denote the set of increasing sequences $i_0 < i_1 < \cdots < i_k$ with $k+1$ elements taken from the set $\{1,2, \ldots ,n\}$.  For each such multi-index
$$I = (i_0, \ldots ,i_k) \in {\cal I}_k, \quad \hbox{we let} \quad G(I) = G(i_0, \ldots ,i_k)$$
denote the subgroup of $G$ generated by $g_{i_0}$, \ldots , $g_{i_k}$.  Given a multi-index $I = (i_0, \ldots ,i_k) \in {\cal I}_k$, we let
$$I^* = (j_0, \ldots ,j_{n-k-2}) \in {\cal I}_{n-k-2}$$
be the complementary multi-index determined by the conditions
\begin{enumerate}
\item $\{ i_0, \ldots ,i_k \} \cap \{ j_0, \ldots ,j_{n-k-2} \} = \emptyset$,
\item $\{ i_0, \ldots ,i_k \} \cup \{ j_0, \ldots ,j_{n-k-2} \} = \{ 1, \ldots , n \}$.
\end{enumerate}
Then $G(I^*)$ is the subgroup of elements $g \in G$ such that
$$g(i) = 1, \qquad \hbox{for $i \in \{i_0,\ldots ,i_k\}$},$$
and we get a direct product of abelian groups
$$G = G(I) \times G(I^*), \quad \hbox{for each} \quad I \in {\cal I}_k.$$
Since the action of $G$ on the polyhedron $P_a$ is given by $g({\bf a}_i) = g(i) {\bf a}_i$, the restriction of this action to $G(I)$ takes the subpolyhedron $P_a^I = P_a^{i_0,\ldots ,i_k}$ to itself, where
\begin{equation} P_a^I = P_a^{i_0,\ldots ,i_k} = \{ p \in P_a : \hbox{ $g(p) = p$ for $g \in G(I^*)$ } \}. \label{E:subpoly} \end{equation}

We claim that the group $G$ acts on cells in dual CW decompositions of $\partial P_a$ and $T^1_{p_0}M$, where $p_0$ is the point in $P_a$ with coordinates
$$y_1 =  \cdots = y_n = 0.$$
The first of these is a triangulation of $\partial P_a$, in which we let
$$\Delta (I) = \Delta (i_0, \ldots ,i_k), \quad \hbox{for} \quad I = (i_0, \ldots ,i_k) \in {\cal I}_k,$$
be the oriented $k$-simplex in $\partial P_a$ spanned by ${\bf a}_{i_0}, \ldots , {\bf a}_{i_k}$.  Moreover, if $g \in G$, we let $g \Delta (I) = g\Delta (i_0, \ldots ,i_k)$ denote the oriented $k$-simplex in principal coordinates spanned by $g(i_0){\bf a}_{i_0}, \ldots , g(i_k){\bf a}_{i_k}$.  Finally, if $\sigma \in S_{k+1}$, the symmetric group of bijections from $\{0,\ldots ,k\}$ to itself, we let
$$\Delta (i_{\sigma (0)}, \ldots ,i_{\sigma (k)}) = (\hbox{sgn } \sigma )\Delta (i_0, \ldots ,i_k), \qquad \hbox{when $i_0 < \cdots < i_k$},$$
and for completeness, we set $\Delta (i_0, \ldots ,i_k) = 0$ when two of the indices $i_0, \ldots ,i_k$ are equal.  The simplex $\Delta (i_0, \ldots ,i_k)$ lies in the $(k+1)$-dimensional coordinate plane spanned by
$$y_i, \quad \hbox{for} \quad i \notin\{i_0, \ldots ,i_k\}.$$
Moreover,
$$g \mapsto g \Delta (I) = g \Delta (i_0, \ldots ,i_k)$$
gives a faithful representation of $G(I)$ onto the $G$-orbit of $\Delta (I)$ in the space of $k$-simplices, which has cardinality $2^k$, while $G(I^*)$ acts trivially on $\Delta (I)$.

As $k$ ranges from $0$ to $n-1$, $I$ ranges over ${\cal I}_k$ and $g$ ranges over $G(I)$, the oriented simplices $\epsilon (g) g \Delta (I)$ give a simplicial decomposition of the boundary of $P_a$ with $2^{k+1} \binom{n}{k+1}$ simplices of dimension $k$.  We can write
\begin{equation}\partial (P_a) = \sum _{g\in G} \epsilon (g) g\Delta(1,\ldots ,n),\label{E:boundary}\end{equation}
regarding it as an expression for the fundamental class of $\partial (P_a)$.

The distinct simplices among $g\Delta(I)$, for $I \in {\cal I}_k$ and $g \in G(I)$, form a basis for a chain complex $\Delta_*(\partial P_a)$ with boundary operator $\partial _\Delta$ defined by
\begin{multline*} \partial _\Delta \Delta (i_0, \ldots ,i_k) \\ = \Delta (i_1, \ldots ,i_k) -  \Delta (i_0, i_2, \ldots ,i_k) + \ldots + (-1)^k \Delta (i_0, \ldots ,i_{k-1}),\end{multline*}
together with the requirement that the boundary operator $\partial _\Delta$ commute with the action of the group $G$.

To each $(n-1)$-simplex $g \Delta (1, \dots  ,n)$, there is an {\em associated asymptotic vector field\/}
$$g(X) = g(1) x_1 e_1 + \cdots + g(n) x_n e_n.$$
Moreover, to each $(n-k)$-simplex there are $2^{k-1}$ associated asymptotic vector fields, corresponding to the $2^{k-1}$ $(n-1)$-simplices which have the given $(n-k)$-simplex in their closures.  

There is a corresponding \lq\lq dual" CW decomposition of the unit sphere $T^1_{p}M$ of unit tangent vectors within $T_{p}M$, for each $p$ in the domain $U$ of our distinguished principal coordinate system.  We first define
$$\Xi(i)(p) = \hbox{Closure of }\{ v \in T_pM :  dy_i(p)(v) \geq | dy_j(p)(v)|, \hbox{for $j \neq i$} \},$$
and its reflection $g_i \Xi(i)(p)$ under $g_i \in G$.  More generally, we define
$$\Xi(I)(p) = \Xi(i_0)(p) \cap \cdots \cap \Xi(i_k)(p), \quad \hbox{for} \quad I = (i_0, \ldots ,i_k) \in {\cal I}_k,$$
and its reflections $g\Xi(I)(p)$ under $g \in G$.  We use the symbol $\Xi(I)$ to represent the family
$$p \mapsto \Xi(I)(p)$$
as $p$ ranges over the domain $U$.  The symmetry group $G$ acts on these families, giving $g\Xi(I)$, for $I \in {\cal I}_k$ and $g \in G$.

Finally, we intersect with the unit tangent space $T^1_pM$ at $p \in M$, and set
$$\square(i) (p) = \Xi (i)(p) \cap T^1_pM, \qquad g_i\square(i) (p) = g_i\Xi(i) (p) \cap T^1_pM.$$
The cells $\square(i)$ and $g_i\square(i)$, for $ 1 \leq i \leq n$, are the $(n-1)$-dimensional cells in a regular CW decomposition $\square_*(T^1_pM)$ of the $(n-1)$-sphere $T^1_pM$, and they correspond to the duals of the $2n$ vertices of the asymptotic polyhedron $P_a$.  The cells of dimension $(n-1)-k$ in this CW decomposition are generated by the cells
$$\square(I)(p) = \square(i_0, \ldots , i_k)(p) =  \Xi(i_0, \ldots , i_k)(p) \cap T^1_pM, \quad \hbox{for $I = (i_0, \ldots ,i_k) \in {\cal I}_k$,}$$
under the action of the symmetry group $G$.  As before, we let
$$\square (i_{\sigma (0)}, \ldots ,i_{\sigma (k)}) = (\hbox{sgn } \sigma ) \square(i_0, \ldots ,i_k), \qquad \hbox{when $i_0 < \cdots < i_k$},$$
and we set $\square (i_0, \ldots ,i_k) = 0$ when two of the indices $i_0, \ldots ,i_k$ are equal.  We regard
$$g \square(I) \quad \hbox{as the cell dual to} \quad g \Delta(I).$$
Once again, the maps $g \mapsto g \square(I)$ give faithful representations of $G(I)$ onto the orbits containing $\square(I)$ in the space of $((n-1)-k)$-cells, which again has cardinaity $2^k$, while $G(I^*)$ acts trivially on $\square(I)$.

If $\square_{n-1-k}(T^1_pM)$ denotes the resulting space of $(n-1-k)$-dimensional cubical cells in the CW decomposition of $T^1_pM$, we can define a boundary map
$$\partial_\Box : \square_{n-1-k}(T^1_pM) \longrightarrow \square_{n-1-k-1}(T^1_pM)$$
by
$$\partial_\Box (\square (i_0, \ldots ,i_k)) = \sum _{i = 1}^n (1 - g_i) \square (i,i_0, \ldots ,i_k),$$
where $\square (i,i_0, \ldots ,i_k) = 0$ if $i \in \{ i_0, \ldots , i_k \}$, together with the requirement that $\partial_\Box$ commute with the action of the group $G$.  This is the boundary operator for a chain complex $\square_*(T^1_pM)$ which gives the cellular homology of $T^1_pM$ as described in \S 2.2 of Hatcher \cite{Hat}.

There is a duality map between the space $\Delta _k(\partial P_a)$ of $k$-simplices in the simplicial subdivision of $\partial P_a$ and the space $\square_{n-1-k}(T^1_pM)$ of $(n-1-k)$-dimensional cubical cells in the CW decomposition of $T^1_pM$, defined by
$$\langle g \square(I) , h\Delta (J) \rangle \\ = \begin{cases} 1, & \hbox{if $I = J \in {\cal I}_k$ and $g = h \in G(I)$}, \cr 0, & \hbox{otherwise.}  \end{cases} $$

\vskip .1 in
\noindent
{\bf Lemma~6.1.} {\sl The boundary maps $\partial _\Box$ and $\partial _\Delta$ are related by the formula
\begin{equation}\langle \partial_\Box (g \square(I)) , h\Delta (J) \rangle = \langle g \square(I) , \partial _\Delta (h\Delta (J)) \rangle, \label{E:dualboundary}\end{equation} 
for $g,h \in G$, $I \in {\cal I}_k$ and $J \in {\cal I}_{k+1}$.}

\vskip .1 in
\noindent
Proof:  We compare the calculations
\begin{multline*}\langle \partial_\Box (\square(i_0, \ldots ,i_k)) , \Delta (i,i_0, \ldots , i_k) \rangle \\ = \langle \partial_\Box (\square(i_0, \ldots ,i_k)) , - g_i\Delta (i,i_0, \ldots , i_k) \rangle = 1 \end{multline*}
and
\begin{multline*} \langle \square(i_0, \ldots ,i_k) , \partial _\Delta (\Delta (i,i_0, \ldots , i_k)) \rangle \\ = \langle \square(i_0, \ldots ,i_k) , - g_i \partial _\Delta (\Delta (i,i_0, \ldots , i_k)) \rangle = 1,\end{multline*}
whenever $i \notin \{ i_0, \ldots, i_k \}$.  Here we have used our conventions on changing signs when reordering indices in $\Delta (i,i_0, \ldots , i_k)$, as well as the fact that $g_i$ changes the orientation on $\square (i,i_1, \ldots , i_k)$ within$\partial_\Box (\square(i_0, \ldots ,i_k))$ or on $\Delta (i,i_0, \ldots , i_k)$ within $\partial _\Delta (\Delta (i,i_0, \ldots , i_k))$ whenever $i \notin \{ i_0, \ldots, i_k \}$.  QED

\vskip .1 in
\noindent
We can use the principal coordinate system to construct a homotopy equivalence from $\partial P_a$ to $T^1_{p_0}M \cong S^{n-1}$ as $a \rightarrow 0$ where $p_0$ is the center of the asymptotic polyhedron.  When we do this, $\square(I)(p_0)$ and $\Delta(I)$ are realized as dual cells in dual CW decompositions of $T^1_{p_0}M$ and $S^{n-1}$.

\section{Establishing the Allendoerfer-Weil formula by integration over the fiber}
\label{S:faces}

To obtain the desired Hazzidakis formulae for the Riemannian polyhedra $P_a$ of (\ref{E:Riemannianpolyhedron}), we need to modify the Gauss-Bonnet formula for manifolds with smooth boundary from \S \ref{S:integrationoverfibers} so that it allows for for abrupt changes in the vector field $V$ ococuring in the geodesic curvature form along simplices of various dimensions in the simplicial decomposition of the boundary $\partial P_a$.

The asymptotic polyhedron $P_a$ is homeomorphic to an $n$-disk with boundary an $(n-1)$-dimensional sphere and lies in the domain $U$ of a principal coordinate system $(y_1, \ldots ,y_n)$, and this defines a trivialization $T^1U = U \times S^{n-1}$.  If the unit-length vector field $V$ along $\partial P_a$ were chosen to be smooth and point outward along every $(n-1)$-simplex of $\partial P_a$ and and to have a single nondegenerate zero of Morse index zero at the center of $P_a$, then it would follow from Remark~5.6 that
$$ \frac{1}{(2\pi)^{n/2}} \int_{P_a} \hbox{Pf}(\Omega ) - \int_{\partial P_a} V^*Te(\nabla ) = 1, \quad \hbox{for $n = \dim M$ is even, and}$$
$$- \int_{\partial P_a} V^*Te(\nabla ) = 1, \quad \hbox{for $n$ odd,} $$
where $Te(\nabla )$ is the transgressed Euler form described before.  But as suggested by Remark~5.5, we would like to choose $V$ along each $(n-1)$-simplex to be the outward-pointing unit-length asymptotic vector field associated to that simplex, in the hope that this will simplify the geodesic curvature form
$$\Pi (V, \nabla) = - V^*Te(\nabla )$$
on that simplex.  But to do this, we need to introduce additional corner terms along the $(n-2)$-dimensional faces which account for the change in $V$ when passing from one $(n-1)$-dimensional face to another, and this process needs to be iterated giving corrections for all the $(n-k)$-dimensional faces of the asymptotic polyhedron, for $1 \leq k \leq n-2$.  When these corrections are made, we get a cycle in the double complex for $\partial P_a \times S^{n-1}$ that was described in the previous section.  The additional terms in this cycle are obtained by integration over the fiber, extending the technique we used in \S \ref{S:integrationoverfibers}.

Our approach is to describe what the resulting cycle in the double complex of \S \ref{S:dualCW} for $T^1U = U \times S^{n-1}$ should be, check that it is indeed a cycle, and show that as $a \rightarrow 0$, it deforms to the fundamental cycle of dimension $n-1$ in the cellular chain complex for the $(n-1)$-sphere over the center of the asymptotic polyhedron $P_a$, which is also homotopic to the cycle representing any smooth outward pointing unit-length vector field which is smooth along $\partial P_a$.

Recall that the cells in the product cellular decomposition of $\partial P_a \times S^{n-1}$ within $T^1U$ are
$$g \Delta (i_0, \ldots ,i_k) \times h \square (j_0, \ldots ,j_l),$$
for arbitrary choice of indices $(i_0, \ldots ,i_k) \in {\cal I}_k$ and $(j_0, \ldots ,j_l) \in {\cal I}_l$, and arbitrary choice of $g, h \in G$.  Here the cells $g\Delta (i_0, \ldots ,i_k)$ give a simplicial decomposition of $\partial P_a$ as described in \S \ref{S:dualCW}, while the cells $h\Box (j_0, \ldots ,j_l)$ provide a cellular decomposition of the fiber $S^{n-1}$.  The boundary is $\partial = \partial _1 + \partial _2$, where
$$\partial _1[g(\Delta(i_0,\ldots ,i_k) \times h\square (j_0, \ldots ,j_l)] = [g\partial_\Delta \Delta(i_0,\ldots ,i_k)] \times h\square (j_0, \ldots ,j_l), $$
and
$$\partial _2[g(\Delta(i_0,\ldots ,i_k) \times h\square (j_0, \ldots ,j_l)] = (-1)^k g\Delta(i_0,\ldots ,i_k) \times [h\partial_\Box \square (j_0, \ldots ,j_l)],$$
the boundaries $\partial _1$ and $\partial _2$ being required to commute with the actions of $G$ on the two factors.

\vskip .1 in
\noindent
{\bf Definition.} The {\em fundamental cycle\/} over $\partial P_a \times S^{n-1}$ is
\begin{equation} z = \sum _{k=0}^{n-1} z_{k,n-1-k}, \label{E:fundcycle} \end{equation}
where
$$z_{k,n-1-k} = (-1)^{k(k-1)/2} \sum _{i_0 <  \cdots < i_k} \sum _{g\in G(i_0, \ldots ,i_k)} g\Delta(i_0,\ldots ,i_k) \times g\square (i_0, \ldots ,i_k),$$
for $0 \leq k \leq n-1$.

\vskip .1 in
\noindent
Note that
$$\dim \Delta(i_0,\ldots ,i_k) = k \quad \hbox{and} \quad \dim \square (i_0, \ldots ,i_k) = n-k-1,$$
so each term in the sum (\ref{E:fundcycle}) has dimension $n-1$.

We define an action $\beta $ of $G$ on the fundamental cycle $z$ by letting $\beta (h,z)$ be obtained by replacing  $g\Delta (I) \times g\square(I)$ by $\epsilon (h)g\Delta (I)\times gh\square(I)$, for $h \in G$, where $\epsilon (h)$ is defined by (\ref{E:epsilong}).  We can think of this as an action of $G$ on the fiber of the cell $g\square (I)\times g\Delta (I)$ over its base $g\Delta (I)$, an action $\beta $ that preserves the Riemannian metric on the fibers of $T^1M$.  Multiplication by $\epsilon (h)$ ensures that this action also preserves the orientation of $T^1M$.

\vskip .1 in
\noindent
{\bf Lemma~7.1.} {\sl The fundamental cycle (\ref{E:fundcycle}) is indeed a cycle in the cellular chain complex of the product $\partial P_a \times S^{n-1}$, that is $\partial z = 0$.  When $a \rightarrow 0$, this cycle approaches a representative of the fundamental class of the fiber $S^{n-1}$ over the center of $P_a$.  Moreover, $\beta (z,h)$ is also a cycle, which approaches the fundamental class of $S^{n-1}$ when $a \rightarrow 0$.}

\vskip .1 in
\noindent
Proof:  We use the notation from \S \ref{S:dualCW} in which ${\cal I}_k$ denotes the set of increasing sequences $(i_0,\ldots ,i_k)$ of length $k+1$ taken from the set $\{ 1, \ldots ,n\}$.  We can then write
$$z_{k,n-1-k} = \frac{(-1)^{k(k-1)/2}}{2^{n-1-k}} \sum _{g\in G} \sum _{I \in {\cal I}_k} g\Delta(I) \times g\square (I)$$
for $0\leq k \leq n-1$.  When differentiating with respect to $\partial _1$, we obtain
\begin{multline*} \partial _1[g\Delta(i_0,\ldots ,i_k) \times \epsilon (g)g\square (i_0, \ldots ,i_k)] \\ = \left[ \sum (-1)^j g\Delta(i_0,\ldots \widehat{i_j}, \ldots,i_k) \right] \times \epsilon (g)g\square (i_0, \ldots ,i_k). \end{multline*}
When $I \in {\cal I}_k$ and $J \in {\cal I}_{k-1}$, we write $J \subseteq I$ if the elements of the sequence $J$ are all included in $I$.  We can then conclude that
\begin{equation} \partial _1(z_{k,n-1-k}) = \frac{(-1)^{k(k-1)/2}}{2^{n-1-k}} \sum _{g\in G} \sum _{I \in {\cal I}_k} \left[ \sum _{J \in {\cal I}_{k-1}, J \subseteq I} g\Delta(J) \times g\square (I)\right]. \label{E:partial1}\end{equation}
On the other hand,
$$\partial _2[g\Delta(i_0,\ldots ,i_{k-1}) \times g\Box (i_0, \ldots ,i_{k-1})] \qquad \qquad \qquad \qquad \qquad \qquad $$
$$\qquad = (-1)^k g\Delta(i_0,\ldots ,i_{k-1}) \times \left[ \sum (1-g_i) \Box (i, i_0, \ldots ,i_{k-1}) \right],$$
from which we conclude that
\begin{multline} \partial _2(z_{k-1,n-k}) = 2 \frac{(-1)^{k + \frac{(k-1)(k-2)}{2}}}{2^{n-k}} \sum _{g\in G} \sum _{I \in {\cal I}_k} \left[ \sum _{J \in {\cal I}_{k-1}, J \subseteq I} g\Delta(J) \times g\square (I)\right] \\
= - \frac{(-1)^{k(k-1)/2}}{2^{n-1-k}} \sum _{g\in G} \sum _{I \in {\cal I}_k} \left[ \sum _{J \in {\cal I}_{k-1}, J \subseteq I} g\Delta(J) \times g\square (I)\right].\label{E:partial2}\end{multline}
Comparison of (\ref{E:partial1}) and (\ref{E:partial2}) yields
$$\partial _1(z_{k,n-k}) + \partial _2(z_{k-1,n-k+1}) = 0, \quad \hbox{for} \quad 0 \leq k \leq n-2.$$
On the other hand, one checks directly that
$$\partial _1 z_{0,n-1} = 0, \quad \partial _2 z_{0,n-1} = 0,$$
so $z$ is a cycle as claimed.

To see that the cycle is homologous (in fact even homotopic) to the fiber when $n$ is even, simply let the parameter $a$ in $P_a$ go to zero.  In the limit, all that remains are the terms
\begin{equation} \sum _{i=1}^n [\Delta(i) \times \square (i) + g_i\Delta(i) \times g_i \square (i)] = \{ 0 \} \times \sum _{i=1}^n (\square (i) + (-1)^nA \square (i)), \label{E:homopmfiber} \end{equation}
where $\{ 0 \}$ denotes the center of the asymptotic polyhedron $P_a$ and $A$ is the antipodal map on $S^{n-1}$, which is orientation-preserving if and only if $n$ is even.  Since $\square (i)$ and $g_i \square (i)$ are antipodal spherical angles at the vertices of the polyhedron which all approach the center point of the polyhedron as $a \rightarrow 0$, this reduces to the cellular description of the fiber over $\{ 0 \}$ as described before.

Finally, the \lq\lq moreover" statement follows from the fact that under the action $\beta $, each element $h \in G$ acts by an orientation-preserving invertible cellular map on the second factor of the product CW complex
$$\Delta _*(\partial P_a) \times \square_*(T^1_pM). \quad \hbox{QED}$$

\vskip .1 in
\noindent
The preceding Lemma allows us to modify the formulae of (\ref{E:GBDeven}) and (\ref{E:GBDodd}) so that they apply to a vector field along $\partial P_a$ which has corners along simplices of dimension $\leq n-2$ in the simplicial decomposition $\Delta _*(\partial P_a)$ of $\partial P_a$.  When $n$ is even we obtain the formula,
\begin{equation} \frac{1}{(2\pi)^{n/2}} \int_{P_a} \hbox{Pf}(\Omega ) - \sum _{k=1}^{n-1} \sum _{I \in {\cal I}_k} \sum _{g\in G(I)} (-1)^{k(k-1)/2} \int_{g\Delta(I) \times g\square (I)} Te(\nabla) - 1. \label{E:gbbdryeven1}\end{equation} 
where the sum of integrals on the right represents the evaluation of transgressed Euler form $Te(\nabla)$ on the fundamental cycle $z$ over $\partial P_a$, which we described at the beginning of this section.

We can integrate each of the terms in (\ref{E:gbbdryeven1}) over the fiber $g\square (I)$, using the slant product as described before.  For $I = (i_0, i_1, \ldots, i_k) \in {\cal I}_k$ and $g \in G$, we can set
$$g\tau (\nabla ,I) = g\tau (\nabla , i_0, \ldots ,i_k) = - Te(\nabla) /g\Xi (i_0, \ldots ,i_k),$$
a differential form of degree $k$.  Note that with the change of sign we have
$$g\tau(\nabla, 1, \ldots, n) = \Pi (\nabla ,V) | g \Delta (1, \ldots, n),$$
which is just the geodesic curvature form restricted to $g \Delta (1, \ldots, n)$, while more generally $g\tau (\nabla ,I)$ is obtained by integration over a cell in the CW decomposition of $T^1M$.  We can rewrite (\ref{E:gbbdryeven1}) as
\begin{equation}\frac{1}{(2\pi)^{n/2}} \int_{P_a} \hbox{Pf}(\Omega ) + \sum _{k=1}^{n} \sum _{I \in {\cal I}_k} \sum _{g\in G(I)} (-1)^{k(k-1)/2} \int_{g\Delta(I)} g\tau (\nabla ,I) = 1. \label{E:gbbdryeven2}\end{equation}
This is the {\em Allendoerfer-Weil formula\/} we need for the proof of the Main Theorem in the even-dimensional case.

There are similar formulae for the case in which $n$ is odd.  We choose a distinguished $(n-1)$-simplex in $\partial P_a$ which we denote by $\Delta _0$.  The collection of all $(n-1)$ simplices in $\partial P_a$ is then $\{\epsilon (g) g\Delta_0:g \in G\}$, the $\epsilon (g)$ being included so that they have the right orientation.  Formula (\ref{E:GBDodd}) then becomes
\begin{equation} \sum _{k=1}^{n} \sum _{I \in {\cal I}_k} \sum _{g\in G(I)} (-1)^{k(k-1)/2} \int_{g\Delta(I) \times g\square (I)} Te(\nabla ) = 0. \label{E:gbbdryodd1}\end{equation} 
As before, we can set
$$g\tau (\nabla ,I) = g\tau (\nabla , i_0, \ldots ,i_k) = \tau (\nabla )/g\Xi (i_0, \ldots ,i_k),$$
for $I = (i_0, i_1, \ldots, i_k) \in {\cal I}_k$, and $(\ref{E:gbbdryodd1})$ becomes
\begin{equation}\sum _{k=1}^{n} \sum _{I \in {\cal I}_k} \sum _{g\in G(I)} (-1)^{k(k-1)/2} \int_{g\Delta(I)} g\tau (\nabla ,I) = 1, \label{E:gbbdryodd2}, \end{equation}
which is {\em Allendoerfer-Weil formula\/} we need for the proof of the Main Theorem in the odd-dimensional case.

\vskip .1 in
\noindent
{\bf Theorem~7.2.} {\sl When $n$ is even, the Gauss-Bonnet formula for the Riemannian $n$-cell $P_a$ is either (\ref{E:gbbdryeven1}) or (\ref{E:gbbdryeven2}).  When $n$ is odd, the Gauss-Bonnet formula for the Riemannian $n$-cell $P_a$ is either (\ref{E:gbbdryodd1}) or (\ref{E:gbbdryodd2}).}

\vskip .1 in
\noindent
To make use of these formulae, we will need to know that the zero-dimensional terms are just solid angles.  To see this, we evaluate the top-dimensional cubical cells at the vertices
$$\Delta (i)^+ = {\bf a}_i \quad \hbox{and} \quad \Delta (i)^- = -{\bf a}_i. $$
of the asymptotic polydedron (\ref{E:Riemannianpolyhedron}).  We obtain the {\em dual angles\/}
\begin{equation} \square(i)^+ = \square(i)({\bf a}_i) \quad \hbox{and} \quad \square(i)^- = g_i\square(i)(-{\bf a}_i) \label{E:solidangles}\end{equation}
to the asymptotic polyhedron at the vertices ${\bf a}_i$ and $-{\bf a}_i$.

\vskip .1 in
\noindent
{\bf Lemma 7.3.}\begin{sl}  The integrals over zero-simplices ${\Delta(i)^+}$ and ${\Delta(i)^-}$ in the simplicial decomposition $\Delta _*(\partial P_a)$ of $\partial P_a$ that appear in the above Allendoerfer-Weil formulae are given by
\begin{equation}  \int_{\square(i)^+} Te(\nabla ) = \frac{\hbox{volume of $\square(i)^+$}}{\hbox{volume of $S^{2m-1}$}}, \quad  \int_{\square(i)^-} Te(\nabla ) = \frac{\hbox{volume of $\square(i)^-$}}{\hbox{volume of $S^{2m-1}$}}, \label{E:sqi}\end{equation}
where $S^{n-1}$ has constant curvature one and hence its volume is given by (\ref{E:volSn-1}). \end{sl}

\vskip .1 in
\noindent
Proof:  In these cases we are integrating over a fiber lying over a single point ${\square(i)^+}$ or ${\square(i)^-}$ for some $i$ with $1 \leq i \leq n$.  We can use polar coordinates as we did earlier when calculating the volume of $S^{n-1}$ and we quickly obtain (\ref{E:sqi}).  QED

\vskip .1 in
\noindent
{\bf Remark 7.4.}  There are similar Allendoerfer-Weil formula for each of the subpolyhedra $P_a^I$ defined by (\ref{E:subpoly}).

\section{Hazzidakis formulae}
\label{S:hazzidakis}

Our goal is to show that after averaging over a suitable family of metric connections in $TM$ most of the terms in the Allendoerfer-Weil formulae presented in the previous section vanish, leaving the Hazzidakis formulae we need for the proof of our Main Theorem.  There are two such Hazzidakis formulae, depending on whether $n$ is even or odd.

\vskip .1 in
\noindent
{\bf Even-dimensional Hazzidakis~Formula.} {\sl Let $f: M^{2m} \rightarrow {\mathbb E}^{4m-1}$ be an isometric immersion, where $M^{2m}$ is a $(2m)$-dimensional Riemannian manifold of constant curvature $-1$.  If $\Omega $ is the curvature of the Levi-Civita connection on $TM$, then the integral of $\hbox{Pf}(\Omega )$ over the asymptotic polyhedron $P_a$ within $M^{2m}$ is given by the formula
\begin{equation} \frac{1}{(2\pi)^m} \int _{P_a} \hbox{Pf}(\Omega ) = 1 - \sum _{i = 1}^n \left[ \frac{\hbox{volume of $\square(i)^+$}}{\hbox{volume of $S^{2m-1}$}} + \frac{\hbox{volume of $\square(i)^-$}}{\hbox{volume of $S^{2m-1}$}} \right]. \label{E:hazzeven}\end{equation}}

\noindent
Since no volume of a solid angle can be larger than the volume of the unit $(n-1)$-sphere, equation (\ref{E:pf=vol}) shows that this Hazzidakis formula gives an explicit bound on the volume of $P_a$.  Moreover, one easily checks using (\ref{E:pf=vol}) that formula (\ref{E:hazzeven}) is a generalization of the Hazzidakis formula for surfaces (\ref{E:hilb12}).

The explicit volume bound given by the Hazzidakis formula implies the even-dimensional generalization of Hilbert\rq s theorem, as stated in the Introduction.

\vskip .1 in
\noindent
{\bf Odd-dimensional Hazzidakis~Formula.} {\sl Let $f: M^{2m+1} \rightarrow {\mathbb E}^{4m+1}$ be an isometric immersion, where $M^{2m+1}$ is a $(2m+1)$-dimensional Riemannian manifold of constant curvature $-1$.  If $\Omega $ is the curvature of the Levi-Civita connection on $TM$, then the integral of $\hbox{Pf}(\Omega )$ over the boundary $\partial P_a$ of the asymptotic polyhedron $P_a$ within $M^{2m+1}$ is given by the formula
\begin{equation} \frac{1}{2} \frac{1}{(2\pi)^m} \int _{\partial P_a} \hbox{Pf}(\Omega ) = 1 - \sum _{i = 1}^n \left[ \frac{\hbox{volume of $\square(i)^+$}}{\hbox{volume of $S^{2m}$}} + \frac{\hbox{volume of $\square(i)^-$}}{\hbox{volume of $S^{2m}$}} \right]. \label{E:hazzodd}\end{equation}}

\vskip .1 in
\noindent
Note that the $\hbox{Pf}(\Omega )$ in (\ref{E:hazzodd}) is the $(2m)$-dimensional Pfaffian calculated with respect to the curvature $\Omega $ of the ambient manifold $M^{2m+1}$.  This gives a bound on the $(2m)$-dimensional volume of $\partial P_a$ because the second-fundamental form terms $(\nabla V)^2/4$ for an appropriate choice of $V$ that would normally occur in the geodesic curvature form give by formula (\ref{E:geodesiccurvforboundaryMmod}) vanish, as we will see in Lemma~9.2, so we can apply equation (\ref{E:pf=vol}) to get a bound on the volume of $\partial P_a$. 

To get a bound on the volume of $P_a$ itself we use the isoperimetric inequality of Theorem 34.2.6 from the book of Burago and Zallgaller \cite{BZ}.  From here, the proof of the Main Theorem proceeds as in the even-dimensional case.

Thus once we have proven the above Hazzidakis formulae we will have proven the Main Theorem.  Remark~5.3 gives a proof when $n = 2$, so we can assume for the rest of the article that $n$ is $\geq 3$.

\section{Applying the Allendoerfer-Weil formulae to metric connections in $TM$}
\label{S:connections}

We can apply the Allendoerfer-Weil formulae of \S \ref{S:faces} to any connection $\nabla $ on the oriented tangent bundle $TM$ which is compatible with the Riemannian metric on $M$, not just the Levi-Civita connection which has vanishing torsion.  Recall that the space of such metric connections is an affine space, that is, if $\nabla $ and $\nabla ^*$ are metric connections on $TM$, so is
$$t \nabla + (1 - t) \nabla ^*, \quad \hbox{for any $t \in {\mathbb R}$.}$$

The principal adapted moving frame $(e_1, \ldots , e_{2n})$ chosen in \S \ref{S:adaptedframe} determines a moving frame $(e_1, \ldots , e_n)$ of orthonormal sections of $TM$ which in turn provides a distinguished flat metric connection $d = \nabla ^0$ with respect to which $e_1, \ldots , e_n$ are parallel.  We can regard this metric connection as the origin in the space of metric connections on $TM$, allowing us to define sums and averages of connections.

An especially important metric connection on $TM$ is the Levi-Civita connection, which we can write as $\nabla = d + \omega $, where $\omega $ denotes the skew-symmetric matrix of Levi-Civita connection forms.

Other metric connections can be constructed by pulling back the standard metric connections in the normal bundle to the developable immersion $h \circ f : M^n \rightarrow {\mathbb H}^{2n}$ by one of the unit-length asymptotic vector fields described in \S \ref{S:adaptedframe}.  Indeed, choice of a distinguished asymptotic vector field
$$X = x_1 e_1 + \cdots + x_n e_n.$$
allows us to construct a vector bundle isomorphism
$$F_X : TM \rightarrow N(h \circ f) = (\hbox{normal bundle of $h \circ f$}), \quad F_X(v) = \sum e_\lambda \Phi _\lambda (X,v),$$
which we use to pull the metric connection in the normal bundle back to a connection $\widehat \nabla $ on $TM$.  Since the second fundamental forms $\Phi _\lambda $ can be diagonalized simultaneously, this connection $\widehat \nabla $ is flat.  Moreover, the distinguished asymptotic vector field $X$ is parallel with respect to this connection, a fact expressed by the equation $\widehat \nabla X = 0$.  We will also write this connection as $d + \phi $, where
$$\phi _{ij} = \omega _{n+i,n+j}, \qquad \hbox{which implies} \qquad \widehat \nabla e_j = \sum e_i \phi _{ij}.$$

We let $Y_1, Y_2, \ldots , Y_n$ denote the coordinate vector fields for our distinguished principal coordinate system $(y_1, y_2, \ldots , y_n)$.  For this coordinate system, the vector fields
$$Y_i = \frac{\partial }{\partial y_i}, \quad \hbox{for $1 \leq i \leq n$,}$$
are tangent to the frame vectors in our orthonormal frame $(e_1, e_2, \ldots , e_n)$.  It then follows from (\ref{eq:hilb28}) and (\ref{eq:hilb29}) that
$$\omega _{ij}(Y_i - Y_j) = - \phi _{ij}(Y_i - Y_j), \quad \hbox{for} \quad 1 \leq i < j \leq n,$$
which in turn gives a relation between the covariant derivatives of the asymptotic vector field $X$ with respect to Levi-Civita and normal connections in the direction $Y_i - Y_j$:
\begin{equation} \nabla _{Y_i - Y_j}X = - \widehat \nabla _{Y_i - Y_j}X = 0, \quad \hbox{for} \quad X = x_1e_1 + \cdots + x_ne_n, \label{E:hilb81}\end{equation}
since $X$ is parallel with respect to $\widehat \nabla $.

Recall that our symmetry group $G$ consists of maps
$$g : \{ 1 , \ldots ,n \} \rightarrow \{ -1,1 \},$$
and acts as isometries on the tangent bundle $TM$.  This group also acts on metric connections on $T$M by taking a connection $\nabla $ to the connection $\nabla ^g$ given by the formula
$$\nabla ^g = g(\nabla ) = g \circ \nabla  \circ g, \quad \hbox{or} \quad g(d +  \omega ) = d + g \omega g.$$
We remind the reader that since the isometry $g$ is usually not the differential of an isometry on the base manifold $M$, the Levi-Civita connection $\nabla $ is usually taken to a metric connection $\nabla ^g$ with non-vanishing torsion.  Observe that if $\nabla $ has curvature matrix $\Omega $, then $\nabla ^g$ has curvature matrix
$$\Omega ^g = g \Omega g.$$

The group $G$ acts simply transitively on the $(n-1)$-simplices in the simplicial decomposition of $\partial P_a$, as well as on the unit-length asymptotic vector fields on $U$.  Moreover, if $X$ is the unit length asymptotic vector associated to the $(n-1)$-simplex $\Delta(1, 2, \ldots , n)$, the element $g \in G$ takes $X$ to the unit-length asymptotic vector field $X^g$ associated to $g\Delta(1, 2, \ldots , n)$.  Note that the connection obtained from the normal connection $\widehat \nabla$ by pulling back via $F_{X^g}$ is just
$$\widehat \nabla ^g = g(\widehat \nabla ) = g \circ \widehat \nabla  \circ g.$$
Thus the collection of unit-length asymptotic vector fields corresponds to a collection $\{ \widehat \nabla ^g : g \in G \}$ of flat connections in the normal bundle to $h \circ f$ which are permuted by the action of $G$.

Finally, the zero connection $\nabla ^0$ is fixed under this action of $G$, that is $(\nabla ^0)^g = \nabla ^0$, for $g \in G$.

\vskip .1 in
\noindent
{\bf Lemma 9.1.}\begin{sl} If $\nabla $ and $\widehat \nabla $ are the Levi-Civita connection and the pullback via $X$ of the normal connection respectively, then when restricted to $\Delta(1, 2, \ldots , n)$, we have
$$ \nabla X = - \widehat \nabla X = 0,$$
and when restricted to $g\Delta(1, 2, \ldots , n)$, we have
$$ \nabla X^g = - \widehat \nabla ^gX^g = 0.$$
\end{sl}

\noindent
Proof: This follows immediately from (\ref{E:hilb81}), the fact that the vector fields $Y_j - Y_j$ generate the tangent space to $\Delta(1, 2, \ldots , n)$, and equivariance under $G$.  QED

\vskip .1 in
\noindent
When $n$ is even, Lemma~9.1 allows us to eliminate terms of dimension $n-1$ in the Allendoerfer-Weil formula, and when $n$ is odd, it will eliminate all terms other than those containing only the curvature of $M$.  Indeed, when $n$ is odd-dimensional, the Gauss-Bonnet integrand is zero while the geodesic curvature form on $g\Delta(1, 2, \ldots , n)$ is given by

\vskip .1 in
\noindent
{\bf Lemma 9.2.}\begin{sl} If $M$ is odd-dimensional, the geodesic curvature form for the Levi-Civita connection $\nabla $ restricted to $g\Delta(1, 2, \ldots , n)$ is $(1/2) \hbox{Pf}(\Omega )$, where $\Omega $ is the curvature on $M$.  \end{sl}

\vskip .1 in
\noindent
Proof: Indeed, in the case of $\Delta(1, 2, \ldots , n)$, we set $V = X$ in (\ref{E:geodesiccurvforboundaryMmod}) and the calculations above show that $(\nabla X)^2/4$ vanishes leaving $(1/2) \hbox{Pf}(\Omega )$ in the formula for the geodesic curvature form.  A similar formula holds for $g\Delta(1, 2, \ldots , n)$ by $G$-equivariance.  QED

\vskip .1 in
\noindent
We next consider the one-dimensional terms:

\vskip .1 in
\noindent
{\bf Lemma 9.3.}\begin{sl} If $I = (i_0, i_1)$ is a multi-index in the index set ${\cal I}_1$ and $\nabla $ and $\widehat \nabla $ are the Levi-Civita connection and the pullback via $X$ of the normal connection respectively, the terms corresponding to this index in the Allendoerfer-Weil formulae (\ref{E:gbbdryeven1}) or (\ref{E:gbbdryodd1}) must vanish; that is,
$$\int_{g\Delta(I) \times g\square (I)} Te(\nabla ) = \int_{g\Delta(I) \times g\square (I)} Te(\widehat \nabla ) = 0, \quad \hbox{for $g \in G$.}$$
Moreover, the same conclusion holds for any connection $\nabla ^h$ which is obtained from $\nabla $ by the action of $h \in G$. 
\end{sl}

\vskip .1 in
\noindent
Proof: Let us consider first the special case in which our multi-index $I = (1,2)$ and $g$ is the identity, which connects the vertex in the positive $y_1$-direction with the vertex in the positive $y_2$-direction.  As in the previous Lemma, we let
$$Y_i = \frac{\partial }{\partial y_i}, \quad \hbox{for $1 \leq i \leq n$.}$$
It follows from the description of the dual CW complexes in \S\ref{S:dualCW} that $\Delta (I)$ lies in the boundary of a family of $(n-1)$-simplices $\{\Delta _\alpha : \alpha \in A \}$ with associated unit-length asymptotic vector fields of the form
$$X_\alpha = Y_1 + Y_2 \pm Y_3 \pm \cdots \pm Y_n, \quad \alpha \in A,$$
where the signs in front of $Y_r$, for $3 \leq r \leq n$ are chosen in $2^{n-2}$ ways.  These are the corners of the $(n-2)$-dimensional cell $\square(I)$ in the cellular decomposition of $T^1M$.  Moreover, $Y_1 - Y_2$ is tangent to $\Delta (I)$, and as in the proof of Lemma~9.1, we find that if $\widehat \nabla ^\alpha $ is the pullback of the normal connection via $X_\alpha $,
$$\nabla _{Y_1-Y_2}X_\alpha = \widehat \nabla ^\alpha_{Y_1-Y_2}X_\alpha = 0,$$
so $\nabla X_\alpha = 0$ when restricted to $\Delta (I)$.  If we choose $V$ to be a unit-length vector in $\square(I)$ which is a constant linear combination of these $X_\alpha $, then it is also the case that $\nabla V =0$ on $\Delta (I)$. These facts imply that the integral of $\nabla V$ over the fiber $\square (I)$ must be zero, and hence this geodesic curvature form simplifies to
\begin{equation} \Pi(\nabla , V) = \frac{(-1)^{[(n-1)/2]}}{\pi^{(n/2)}} \int _0^\infty e^{-t^2} dt \hbox{Tr}^M_s[V \exp ({\mathcal R}^M )].\label{E:fiberintnew} \end{equation}
This geodesic curvature form must vanish when integrated over $\Delta(I) \times \square (I)$, since the base $\Delta (I)$ is one-dimensional and ${\mathcal R}^M$ contains the factors $\theta _i \wedge \theta_j$ which are two-dimensional in the base.

The more general case in which $I$ is any element of ${\cal I}_1$ and $g$ is any element of $G$ follows by a possible renumbering of the prinicipal coordinates and $G$-equivariance.  QED

\vskip .1 in
\noindent
Unfortunately the technique we have just applied doesn't quite work for terms of dimension $k$ when $1 < k < n-1$.  However, the argument of Lemma~9.3 can be modified to provide a considerable simplification of the geodesic curvature form (\ref{E:fiberint}) when restricted to $\partial P_a$:

\vskip .1 in
\noindent
{\bf Lemma 9.4.}\begin{sl} When restricted to a cell $g\Delta(I) \times g\square (I)$ of $\partial P_a$ where $I \in {\mathcal I}_k$, the geodesic curvature form (\ref{E:fiberint}) simplifies to
\begin{multline}  \Pi(\nabla , V)|g\Delta(I) \times g\square (I) \\ = \frac{(-1)^{[(n-1)/2]}}{\pi^{(n/2)}} \int _0^\infty e^{-t^2} dt \hbox{Tr}^{g\Delta(I)}_s[V \exp ({\mathcal R}^M )], \label{E:lemma941a}\end{multline}
where $\hbox{Tr}^{g\Delta(I)}_s$ is the supertrace on $g\Delta (I)$, a smooth manifold with a boundary that has corners.  Moreover, the right-hand side simplifies further, yielding
\begin{equation}  \Pi(\nabla , V)|g\Delta(I) \times g\square (I) \\ = c_k \frac{\hbox{volume of $\square(i)^+$}}{\hbox{volume of $S^{2m}$}} \int_{g\Delta(I)} \hbox{Tr}^{g\Delta(I)}_s[ \exp ({\mathcal R}^M )],\label{E:lemma941b} \end{equation}
where $c_k$ is a constant depending only on $k \in {\mathbb N}$ which vanishes when $k$ is odd.
\end{sl}

\vskip .1 in
\noindent
Proof:  The proof is quite similar to the proof of Lemma~9.3.  We can reduce to the special case in which 
$$I = (1,2, \ldots , k+1) \in {\mathcal I}_k \quad \hbox{and $g$ is the identity,}$$
by renumbering the principal coordinates and using $G$-equivariance.  The vertices of $\square (I)$ are then
$$X_\alpha = Y_1 + Y_2 + \cdots +Y_{k+1} \pm Y_{k+2} \pm \cdots \pm Y_n, \quad \alpha \in A,$$
where the signs in front of $Y_r$, for $3 \leq r \leq n$ are chosen in $2^{n-k-1}$ ways.  These are the corners of the $(n-k-2)$-dimensional cell $\square(I)$ in the cellular decomposition of $T^1M$.  Moreover,
$$Y_i - Y_j, \quad \hbox{for} \quad 1 \leq i < j \leq k+1,$$
is tangent to $\Delta (I)$, and hence
$$\nabla _{Y_i-Y_j}X_\alpha = 0, \quad \hbox{for} \quad \alpha \in A.$$
Thus the argument for the preceding lemma shows that the geodesic curvature form reduces to (\ref{E:fiberintnew}).

This time, however, the curvature factor does not immediately vanish.  However, $\exp ({\mathcal R}^M )$ depends only on the coordinates in the base $g \Delta (I)$ not the fiber $g\square (I)$, while the integrations of $V$ over $g \square (I)$ can be evaluating in exactly the same way as in the proof of Lemma~7.3.   This implies that (\ref{E:lemma941a}) reduces to (\ref{E:lemma941b}). QED

\vskip .1 in
\noindent
In conclusion, there are curvature terms in every even dimension that might not be zero when using the Levi-Civita connection $\nabla $ in the Allendoerfer-Weil formulae.

\section{Signed averages and the proof of the Main Theorem}
\label{S:proof}

The final step in our proof involves taking the signed average of the Allendoerfer-Weil formulae of \S \ref{S:faces} over the family $\{ \nabla ^g : g \in G \}$ of metric connections related to the Levi-Civita connection by the action of the symmetry group $G$ described at the end of \S \ref{S:asymptotic}.  In taking the signed average, we simply multiply the Allendoerfer-Weil formula for $\nabla ^g$ by
$$\epsilon (g) = \left\{\begin{array}{ll} 1 & \mbox{if $g$ is orientation-preserving,} \\ -1 & \mbox{if $g$ is orientation-reversing,} \end{array} \right. $$
and then take the average of the results.

We first note that after multiplying by $\epsilon (g)$, the terms of dimensions zero, one, $n-1$ or $n$ remain the same.  Indeed, the $n$-dimensional term is just the integral of the Gauss-Bonnet integrand over $P_a$ and this is preserved since the Euler class is independent of positively oriented metric connection, while multiplication by $\epsilon (g)$ provides the necessary correction for negatively oriented metric connections.  For the $(n-1)$-dimensional terms we note that $G$ acts simply transitively on the $(n-1)$-dimensional faces $g \Delta (1, \ldots , n)$ so the action $\nabla \mapsto \nabla ^g$ simply permutes the terms corresponding to these faces, and multiplication by $\epsilon (g)$ once again provides the needed correction when $g$ is oreintation-reversing.

Similarly, one easily checks that the action of $G$ on $\square (i)$ either takes  $\square (i)$ to itself or to $A(\square (i))$, where $A = g_1 \cdots g_n$ is the antipodal map on $S^{n-1}$.  This antipodal map is either orientation-preserving or orientation-reversing depending on whether $n$ is even or odd, and multiplying by $\epsilon (A)$ reverses signs in the orientation-reversing case.  Thus the zero-dimensional terms are preserved.  Moreover, the one-dimensional terms vanish, as verified by the proof of Lemma~9.3.

Finally, we claim that taking the average over our family of metric connections eliminates all the terms of intermediate dimension $k$, for $1 < k< n-1$.

To see this, we let $I = (i_0, i_1, \ldots, i_k)$ be a multi-index in the index set ${\cal I}_k$ as described at the beginning of \S \ref{S:dualCW}, and divide $G$ into a direct product $G = G(I) \times G'(I)$, with $G'(I)$ fixing the axes $(y_{i_0}, \ldots , y_{i_k})$ in the principal coordinate system $(y_1, \ldots , y_n)$.  We then consider the effect of replacing $\nabla $ by $\nabla ^g$ for $g \in G$ in the formula (\ref{E:lemma941b}).  The sign change in this term is the same as the sign change in the supertrace on $\Delta (I)$, which is given by
$$\epsilon_I(g) = \left\{\begin{array}{ll} 1 & \mbox{if $g|\Delta (I)$ is orientation-preserving,} \\ -1 & \mbox{if $g|\Delta (I)$ is orientation-reversing.} \end{array} \right. $$
But $\epsilon_I(g)$ differs from $\epsilon (g)$ in exactly half the cases when $1 < k < n-1$, so when we average using $\epsilon (g)$ instead of $\epsilon_I(g)$, we get zero.

In other words, when we take the signed average of the Allendoerfer -Weil fromulae over the connections in the family $\{ \nabla ^g : g \in G \}$, the highest and lowest dimensional terms remain exactly the same as for the Levi-Civita connection, but the terms of intermediate dimension vanish.  Thus we obtain the Hazzidakis formulae of \S \ref{S:hazzidakis}, which proves the Main Theorem as we have seen before.

\end{document}